\documentclass[12pt]{article}

\usepackage{amsmath,amssymb,amsthm,amscd,a4wide,upref}

\usepackage[enableskew]{youngtab}

\usepackage{hyperref}
\hypersetup{colorlinks,citecolor=blue,filecolor=black,linkcolor=blue,urlcolor=blue}
\usepackage{enumerate}
\usepackage{mathtools}

\usepackage{lmodern}     % set math font to Latin modern math
\usepackage[T1]{fontenc}
\begin{document}

% % % % % % % % % % % % % %
%\newcommand{\tma}{\textcolor{magenta}}
%\newcommand{\tre}{\textcolor{red}}
% % % % % % % % % % % % % %

\newcommand{\ad}{{\rm ad}}
\newcommand{\cri}{{\rm cri}}
\newcommand{\row}{{\rm row}}
\newcommand{\col}{{\rm col}}
\newcommand{\End}{{\rm{End}\ts}}
\newcommand{\Rep}{{\rm{Rep}\ts}}
\newcommand{\Hom}{{\rm{Hom}}}
\newcommand{\Mat}{{\rm{Mat}}}
\newcommand{\ch}{{\rm{ch}\ts}}
\newcommand{\chara}{{\rm{char}\ts}}
\newcommand{\diag}{{\rm diag}}
\newcommand{\non}{\nonumber}
\newcommand{\wt}{\widetilde}
\newcommand{\wh}{\widehat}
\newcommand{\ol}{\overline}
\newcommand{\ot}{\otimes}
\newcommand{\la}{\lambda}
\newcommand{\La}{\Lambda}
\newcommand{\De}{\Delta}
\newcommand{\al}{\alpha}
\newcommand{\be}{\beta}
\newcommand{\ga}{\gamma}
\newcommand{\Ga}{\Gamma}
\newcommand{\ep}{\epsilon}
\newcommand{\ka}{\kappa}
\newcommand{\vk}{\varkappa}
\newcommand{\si}{\sigma}
\newcommand{\vs}{\varsigma}
\newcommand{\vp}{\varphi}
\newcommand{\de}{\delta}
\newcommand{\ze}{\zeta}
\newcommand{\om}{\omega}
\newcommand{\Om}{\Omega}
\newcommand{\ee}{\epsilon^{}}
\newcommand{\su}{s^{}}
\newcommand{\hra}{\hookrightarrow}
\newcommand{\ve}{\varepsilon}
\newcommand{\ts}{\,}
\newcommand{\vac}{\mathbf{1}}
\newcommand{\di}{\partial}
\newcommand{\qin}{q^{-1}}
\newcommand{\tss}{\hspace{1pt}}
\newcommand{\Sr}{ {\rm S}}
\newcommand{\U}{ {\rm U}}
\newcommand{\BL}{ {\overline L}}
\newcommand{\BE}{ {\overline E}}
\newcommand{\BP}{ {\overline P}}
\newcommand{\AAb}{\mathbb{A}\tss}
\newcommand{\CC}{\mathbb{C}\tss}
\newcommand{\KK}{\mathbb{K}\tss}
\newcommand{\QQ}{\mathbb{Q}\tss}
\newcommand{\SSb}{\mathbb{S}\tss}
\newcommand{\TT}{\mathbb{T}\tss}
\newcommand{\ZZ}{\mathbb{Z}\tss}
\newcommand{\DY}{ {\rm DY}}
\newcommand{\X}{ {\rm X}}
\newcommand{\Y}{ {\rm Y}}
\newcommand{\Z}{{\rm Z}}
\newcommand{\Ac}{\mathcal{A}}
\newcommand{\Lc}{\mathcal{L}}
\newcommand{\Mc}{\mathcal{M}}
\newcommand{\Pc}{\mathcal{P}}
\newcommand{\Qc}{\mathcal{Q}}
\newcommand{\Rc}{\mathcal{R}}
\newcommand{\Sc}{\mathcal{S}}
\newcommand{\Tc}{\mathcal{T}}
\newcommand{\Bc}{\mathcal{B}}
\newcommand{\Ec}{\mathcal{E}}
\newcommand{\Fc}{\mathcal{F}}
\newcommand{\Gc}{\mathcal{G}}
\newcommand{\Hc}{\mathcal{H}}
\newcommand{\Uc}{\mathcal{U}}
\newcommand{\Vc}{\mathcal{V}}
\newcommand{\Wc}{\mathcal{W}}
\newcommand{\Yc}{\mathcal{Y}}
\newcommand{\Ar}{{\rm A}}
\newcommand{\Br}{{\rm B}}
\newcommand{\Ir}{{\rm I}}
\newcommand{\Fr}{{\rm F}}
\newcommand{\Jr}{{\rm J}}
\newcommand{\Or}{{\rm O}}
\newcommand{\GL}{{\rm GL}}
\newcommand{\Spr}{{\rm Sp}}
\newcommand{\Rr}{{\rm R}}
\newcommand{\Zr}{{\rm Z}}
\newcommand{\gl}{\mathfrak{gl}}
\newcommand{\middd}{{\rm mid}}
\newcommand{\ev}{{\rm ev}}
\newcommand{\Pf}{{\rm Pf}}
\newcommand{\Norm}{{\rm Norm\tss}}
\newcommand{\oa}{\mathfrak{o}}
\newcommand{\spa}{\mathfrak{sp}}
\newcommand{\osp}{\mathfrak{osp}}
\newcommand{\f}{\mathfrak{f}}
\newcommand{\g}{\mathfrak{g}}
\newcommand{\h}{\mathfrak h}
\newcommand{\n}{\mathfrak n}
\newcommand{\z}{\mathfrak{z}}
\newcommand{\Zgot}{\mathfrak{Z}}
\newcommand{\p}{\mathfrak{p}}
\newcommand{\sll}{\mathfrak{sl}}
\newcommand{\agot}{\mathfrak{a}}
\newcommand{\qdet}{ {\rm qdet}\ts}
\newcommand{\Ber}{ {\rm Ber}\ts}
\newcommand{\HC}{ {\mathcal HC}}
\newcommand{\cdet}{{\rm cdet}}
\newcommand{\rdet}{{\rm rdet}}
\newcommand{\tr}{ {\rm tr}}
\newcommand{\gr}{ {\rm gr}\ts}
\newcommand{\str}{ {\rm str}}
\newcommand{\loc}{{\rm loc}}
\newcommand{\Gr}{{\rm G}}
\newcommand{\sgn}{ {\rm sgn}\ts}
\newcommand{\sign}{{\rm sgn}}
\newcommand{\ba}{\bar{a}}
\newcommand{\bb}{\bar{b}}
\newcommand{\bi}{\bar{\imath}}
\newcommand{\bj}{\bar{\jmath}}
\newcommand{\bk}{\bar{k}}
\newcommand{\bl}{\bar{l}}
\newcommand{\hb}{\mathbf{h}}
\newcommand{\Sym}{\mathfrak S}
\newcommand{\fand}{\quad\text{and}\quad}
\newcommand{\Fand}{\qquad\text{and}\qquad}
\newcommand{\For}{\qquad\text{or}\qquad}
\newcommand{\for}{\quad\text{or}\quad}
\newcommand{\grpr}{{\rm gr}^{\tss\prime}\ts}
\newcommand{\degpr}{{\rm deg}^{\tss\prime}\tss}

\renewcommand{\theequation}{\arabic{section}.\arabic{equation}}

\numberwithin{equation}{section}

\newtheorem{thm}{Theorem}[section]
\newtheorem{lem}[thm]{Lemma}
\newtheorem{prop}[thm]{Proposition}
\newtheorem{cor}[thm]{Corollary}
\newtheorem{conj}[thm]{Conjecture}
\newtheorem*{mthm}{Main Theorem}
\newtheorem*{mthma}{Theorem A}
\newtheorem*{mthmb}{Theorem B}
\newtheorem*{mthmc}{Theorem C}
\newtheorem*{mthmd}{Theorem D}

\theoremstyle{definition}
\newtheorem{defin}[thm]{Definition}

\theoremstyle{remark}
\newtheorem{remark}[thm]{Remark}
\newtheorem{example}[thm]{Example}
\newtheorem{examples}[thm]{Examples}

\newcommand{\bth}{\begin{thm}}
\renewcommand{\eth}{\end{thm}}
\newcommand{\bpr}{\begin{prop}}
\newcommand{\epr}{\end{prop}}
\newcommand{\ble}{\begin{lem}}
\newcommand{\ele}{\end{lem}}
\newcommand{\bco}{\begin{cor}}
\newcommand{\eco}{\end{cor}}
\newcommand{\bde}{\begin{defin}}
\newcommand{\ede}{\end{defin}}
\newcommand{\bex}{\begin{example}}
\newcommand{\eex}{\end{example}}
\newcommand{\bes}{\begin{examples}}
\newcommand{\ees}{\end{examples}}
\newcommand{\bre}{\begin{remark}}
\newcommand{\ere}{\end{remark}}
\newcommand{\bcj}{\begin{conj}}
\newcommand{\ecj}{\end{conj}}

\newcommand{\bal}{\begin{aligned}}
\newcommand{\eal}{\end{aligned}}
\newcommand{\beq}{\begin{equation}}
\newcommand{\eeq}{\end{equation}}
\newcommand{\ben}{\begin{equation*}}
\newcommand{\een}{\end{equation*}}

\newcommand{\bpf}{\begin{proof}}
\newcommand{\epf}{\end{proof}}

\def\beql#1{\begin{equation}\label{#1}}

%%% Slaven's notation

%\newcommand{\End}{\mathop{\mathrm{End}}}
%\newcommand{\Hom}{\mathop{\mathrm{Hom}}}
%\newcommand{\vac}{\mathop{\mathrm{\boldsymbol{1}}}}
\newcommand{\Res}{\mathop{\mathrm{Res}}}

\title{\Large\bf Center at the critical level for centralizers in type $A$}

\author{A. I. Molev}

\date{} % Start April 2019
\maketitle

%\vspace{4 mm}

\begin{abstract}
We consider the affine vertex algebra at the critical level
associated with the centralizer of a nilpotent element in the Lie algebra $\gl_N$.
Due to a recent result of Arakawa and Premet, the center of this vertex algebra
is an algebra of polynomials. We construct a family of free generators of the center
in an explicit form. As a corollary, we obtain generators of the
corresponding quantum shift of argument subalgebras
and recover free generators of the center of
the universal enveloping algebra of the centralizer
produced earlier by Brown and Brundan.

\end{abstract}

%\vspace{5 mm}
%%%
%%%{\it Key words:}
%%%

%\newpage

%\tableofcontents
%
%\newpage

\section{Introduction}
\label{sec:int}

Let $\agot$ be a finite-dimensional
Lie algebra over $\CC$ equipped with
an invariant symmetric bilinear form. Consider the corresponding
affine Kac--Moody algebra $\wh\agot$ defined as a central extension of the Lie
algebra of Laurent polynomials $\agot\tss[t,t^{-1}]$; see \eqref{km} below.
The vacuum module over $\wh\agot$ is a vertex algebra whose
center is a commutative associative algebra.
In the case of a simple Lie algebra $\agot$ the structure of the center $\z(\wh\agot)$
at the critical level
was described by a celebrated theorem of Feigin and Frenkel~\cite{ff:ak} (see also~\cite{f:lc})
which states that $\z(\wh\agot)$ is an algebra of polynomials in infinitely many
variables. This fact can be regarded as an affine version
or {\em chiralization} of the well-known description of the center of the universal enveloping
algebra $\U(\agot)$ as a polynomial algebra.

On the other hand, in the case where $\agot=\g^e$ is the centralizer of
a nilpotent element $e$ in a simple Lie algebra $\g$,
the invariant algebra $\Sr(\agot)^{\agot}$ admits free families of generators
under certain additional conditions, due to
the pioneering work of Panyushev, Premet and Yakimova~\cite{ppy:si}.
In particular, these conditions hold for all nilpotents in types $A$ and $C$ which
thus confirms {\em Premet's conjecture} in those cases.
Explicit generators of both the invariant algebra and the center
of $\U(\agot)$ in type $A$ were produced
by Brown and Brundan~\cite{bb:ei} with the use of the shifted Yangians.
It was shown by Yakimova~\cite{o:sp} that the generators of $\Sr(\agot)^{\agot}$
produced in \cite{bb:ei} coincide with those previously conjectured in \cite{ppy:si}.

A recent work of Arakawa and Premet~\cite{ap:qm} provides a chiralization
of these results for the centralizer $\agot$
which shows that the Feigin--Frenkel theorem extends to
the affine vertex algebra at the critical level associated with the
Kac--Moody algebra $\wh\agot$.
As a consequence, they showed
the existence of regular quantum shift of argument subalgebras and proved that they
are free polynomial algebras. Moreover, explicit formulas for generators of $\z(\wh\agot)$
were produced in \cite{ap:qm} in the case where $\agot$ is the centralizer of
a minimal nilpotent in $\gl_N$.

Our goal in this paper is to give explicit formulas for generators of $\z(\wh\agot)$
for an arbitrary nilpotent element $e\in \gl_N$ (Theorem~\ref{thm:suga}).
As a corollary, we get generators
of the quantum shift of argument subalgebras $\Ac_{\chi}\in\U(\agot)$
associated with an element $\chi\in\agot^*$ thus providing
an explicit solution of {\em Vinberg's quantization problem}
\cite{v:sc} for centralizers (Corollary~\ref{cor:genamu}).
It is based on the results of \cite{ap:qm} showing that if the element
$\chi$ is regular, then the classical limit $\gr\Ac_{\chi}$ coincides with the
{\em Mishchenko--Fomenko subalgebra} $\overline\Ac_{\chi}\subset\Sr(\agot)$; see \cite{mf:ee}.
We conjectured in \cite[Conjecture~5.8]{my:qn}
that the subalgebra $\Ac_{\chi}$ can also be obtained with the use of the symmetrization map.
Note that the subalgebra $\Ac_{\tss 0}$ is the center of $\U(\agot)$. In this case
our generators essentially coincide with those found in~\cite{bb:ei}; see Corollary~\ref{cor:gencent}
below.

By taking $e=0$ in Theorem~\ref{thm:suga} we recover
the formulas for the generators of the algebra
$\z(\wh\gl_n)$ found in
\cite{cm:ho} and \cite{ct:qs} which we also interpret with the use
of an operator $\De$ preserving $\z(\wh\gl_n)$ (Corollary~\ref{cor:ff}).
For more details on this particular case and
extensions to the other classical types see \cite{m:so}.
A different way to produce generators of the Feigin--Frenkel center $\z(\wh\g)$ was developed
in a recent work by Yakimova~\cite{o:sf}.

\section{Segal--Sugawara vectors}
\label{sec:ssv}

Using the notation of \cite{bb:ei},
suppose that $e\in\g=\gl_N$ is a nilpotent matrix with Jordan blocks of sizes
$\la_1,\dots,\la_n$, where
$\la_1\leqslant\dots\leqslant \la_n$ and $\la_1+\dots+\la_n=N$.
Consider the corresponding {\em pyramid} which is
a left-justified array of rows of unit boxes such that the top row
contains $\la_1$ boxes, the next row contains $\la_2$ boxes, etc.
The {\em row-tableau} is obtained by writing the numbers $1,\dots,N$ into the boxes
of the pyramid consecutively by rows from left to right.
For instance, the row-tableau
\ben
\young(12,345,6789)
\een
corresponds to the Jordan blocks of sizes $2,3,4$ and $N=9$. Let
$\row(a)$ and $\col(a)$ denote the row and column number of the box containing
the entry $a$.

Let $e_{a\tss b}$ with $a,b=1,\dots,N$ be the standard basis elements of $\g$.
For any $1\leqslant i,j\leqslant n$ and $\la_j-\min(\la_i,\la_j)\leqslant r<\la_j$
set
\ben
E_{ij}^{(r)}=\sum_{\underset{\scriptstyle\col(b)-\col(a)=r}{\row(a)=i,\ \row(b)=j}} e_{a\tss b},
\een
summed over $a,b\in\{1,\dots,N\}$. The elements $E_{ij}^{(r)}$ form a basis
of the Lie algebra $\agot=\g^e$. The commutation relations for
the basis vectors have the form
\beql{comm}
\big[E^{(r)}_{ij},E^{(s)}_{kl}\big]=\de_{kj}\ts E^{(r+s)}_{i\tss l}-\de_{i\tss l}\ts E^{(r+s)}_{kj},
\eeq
assuming that $E^{(r)}_{ij}=0$ for $r\geqslant \la_j$.

Note that in the particular case of a rectangular pyramid $\la_1=\dots=\la_n=p$
the Lie algebra $\agot$ is isomorphic to the truncated polynomial current algebra
$\gl_n[v]/(v^p=0)$ (also known as the {\em Takiff algebra}).
The isomorphism is given by
\ben
E^{(r)}_{ij}\mapsto e_{ij}\ts v^r,\qquad r=0,\dots,p-1,\qquad 1\leqslant i,j\leqslant n.
\een

Equip the Lie algebra $\agot$ with
the invariant symmetric bilinear form $\langle\ts\ts,\ts\rangle$
defined as the normalized Killing form at the critical level
on the $0$-th component of $\agot$; see \cite[Sec.~5]{ap:qm}.
Explicitly, it is given by the formulas
\ben
\big\langle E_{ii}^{(0)},E_{jj}^{(0)}\big\rangle=
\min(\la_i,\la_j)-\de_{ij}\big(\la_1+\dots+\la_{i-1}+(n-i+1)\la_i\big),
\een
and if $\la_i=\la_j$ for some $i\ne j$ then
\ben
\big\langle E_{ij}^{(0)},E_{ji}^{(0)}\big\rangle=-\big(\la_1+\dots+\la_{i-1}+(n-i+1)\la_i\big),
\een
whereas all remaining values of the form on the basis vectors are zero. Note that the sum
in the brackets equals the number of boxes in the first $\la_i$ columns of the pyramid.

The corresponding affine Kac--Moody algebra $\wh\agot$
is the central
extension
\beql{km}
\wh\agot=\agot\tss[t,t^{-1}]\oplus\CC \vac,
\eeq
where $\agot[t,t^{-1}]$ is the Lie algebra of Laurent
polynomials in $t$ with coefficients in $\agot$. For any $r\in\ZZ$ and $X\in\g$
we will write $X[r]=X\ts t^r$. The commutation relations of the Lie algebra $\wh\agot$
have the form
\ben
\big[X[r],Y[s]\big]=[X,Y][r+s]+r\ts\de_{r,-s}\langle X,Y\rangle\ts \vac,
\qquad X, Y\in\agot,
\een
and the element $\vac$ is central in $\wh\agot$.
The {\em vacuum module at the critical level}
over $\wh\agot$
is the quotient
\ben
V(\agot)=\U(\wh\agot)/\Ir,
\een
where $\Ir$ is the left ideal of $\U(\wh\agot)$ generated by $\agot[t]$
and the element $\vac-1$.
By the Poincar\'e--Birkhoff--Witt theorem,
the vacuum module is isomorphic to the universal enveloping algebra
$\U\big(t^{-1}\agot[t^{-1}]\big)$, as a vector space.
This vector space is equipped with a vertex algebra structure; see \cite{fb:va}, \cite{k:va}.
Denote by $\z(\wh\agot)$ the {\em center} of this vertex algebra
which is defined as the subspace
\ben
\z(\wh\agot)=\{v\in V(\agot)\ |\ \agot[t]\tss v=0\}.
\een
It follows from the axioms of vertex algebra that $\z(\wh\agot)$
is a unital commutative associative algebra. It can be regarded as a commutative
subalgebra of $\U\big(t^{-1}\agot[t^{-1}]\big)$.
This subalgebra is invariant with respect to the
{\em translation operator}
$T$ which is
the derivation of the algebra $\U\big(t^{-1}\agot[t^{-1}]\big)$
whose action on the generators is given by
\ben
T:X[r]\mapsto -r\tss X[r-1],\qquad X\in\agot, \quad r<0.
\een

Any element of $\z(\wh\agot)$
is called a {\em Segal--Sugawara vector\/}. By \cite[Thm~1.4]{ap:qm},
there exists a {\em complete set of Segal--Sugawara vectors} $S_1,\dots,S_N$, which means that
all translations $T^r S_l$ with $r\geqslant 0$ and $l=1,\dots,N$
are algebraically independent and
any element of $\z(\wh\agot)$ can be written as a polynomial
in the shifted vectors; that is,
\ben
\z(\wh\agot)=\CC[T^{\tss r}S_l\ |\ l=1,\dots,N,\ \ r\geqslant 0].
\een
In the case $e=0$ this reduces to the Feigin--Frenkel theorem in type $A$~\cite{ff:ak, f:lc}.

To construct a complete set of Segal--Sugawara vectors,
for all $i,j\in\{1,\dots,n\}$ introduce polynomials in a variable $u$
with coefficients in this algebra by
\ben
E_{ij}(u)=\begin{cases}E^{(0)}_{ij}[-1]+\dots+E^{(\la_j-1)}_{ij}[-1]\ts u^{\la_j-1}
&\text{if}\quad i\geqslant j,\\[0.4em]
E^{(\la_j-\la_i)}_{ij}[-1]\tss u^{\la_j-\la_i}+\dots+E^{(\la_j-1)}_{ij}[-1]\ts u^{\la_j-1}
&\text{if}\quad i< j.
\end{cases}
\een
For another variable $x$
calculate the {\em column-determinant}
\begin{align}
\label{cdete}
\cdet\left[\begin{matrix}
x+\la_1\tss T+E_{11}(u)&E_{12}(u)&\dots&E_{1n}(u)\\[0.2em]
E_{21}(u)&x+\la_2\tss T+E_{22}(u)&\dots&E_{2n}(u)\\
\vdots&\vdots& \ddots&\vdots     \\
E_{n1}(u)&E_{n2}(u)&\dots&x+\la_n\tss T+E_{n\tss n}(u)
                \end{matrix}\right]&\\[1.2em]
                \non
{}=\sum_{\si\in\Sym_n} \sgn\si\cdot \big(\de_{\si(1)\tss 1}\ts(x+\la_1\tss T)+E_{\si(1)\tss 1}(u)\big)\dots
\big(\de_{\si(n)\tss n}\ts(x+{}&\la_n\tss T)+E_{\si(n)\tss n}(u)\big)
\end{align}
and write it as a polynomial in $x$ with coefficients in $V(\agot)[u]$,
\beql{polt}
x^n+\phi_1(u)\tss x^{n-1}+\dots+\phi_n(u),\qquad \phi_k(u)=\sum_r\phi^{(r)}_k\ts u^r.
\eeq

\bth\label{thm:suga}
The coefficients $\phi^{(r)}_k$ with $k=1,\dots,n$ and
\beql{conda}
\la_{n-k+2}+\dots+\la_n< r+k\leqslant\la_{n-k+1}+\dots+\la_n
\eeq
belong to the center $\z(\wh\agot)$
of the vertex algebra $V(\agot)$.
Moreover, they form a complete set of
Segal--Sugawara vectors for the Lie algebra $\agot$.
\eth

\bes\label{ex:mn}
In the principal nilpotent case with $n=1$ and $e=e_{1\tss 2}+\dots+e_{N-1\ts N}$ we have
\ben
\phi^{(r)}_1=E^{(r)}_{11}[-1]=e^{}_{1\ts r+1}[-1]+\dots+e^{}_{N-r\ts N}[-1],\qquad r=0,\dots,N-1.
\een

For $n=2$ we have
\ben
\cdet\Bigg[\begin{matrix}
x+\la_1\tss T+E_{11}(u)&E_{12}(u)\\[0.4em]
E_{21}(u)&x+\la_2\tss T+E_{22}(u)
                \end{matrix}\Bigg]=
x^2+\phi_1(u)\tss x+\phi_2(u)
\een
with
\ben
\bal
\phi_1(u)&=E_{11}(u)+E_{22}(u),\\[0.4em]
\phi_2(u)&=E_{11}(u)\tss E_{22}(u)-E_{21}(u)\tss E_{12}(u)+\la_1\ts E^{\ts\prime}_{22}(u),
\eal
\een
where $X'(u)=T\tss X(u)$. Therefore,
a complete set of Segal--Sugawara vectors for $\agot$ is given by
\ben
\phi^{(r)}_1=E^{(r)}_{11}[-1]+E^{(r)}_{22}[-1]\quad \text{with}\quad r=0,1,\dots,\la_2-1,
\een
and
\ben
\phi^{(r)}_2=\sum_{a+b=r}\ts\ts\begin{vmatrix}
E^{(a)}_{11}[-1]&E^{(b)}_{12}[-1]\\[0.4em]
E^{(a)}_{21}[-1]&E^{(b)}_{22}[-1]
                \end{vmatrix}+\la_1\ts E^{(r)}_{22}[-2]
\een
with \ \ $r=\la_2-1,\dots,\la_1+\la_2-2$, where the vertical lines indicate
the column-determinant.

The minimal nilpotent case $e=e_{n\ts n+1}\in\gl_{n+1}$
corresponds to the pyramid with the $n$ rows $1,\dots,1,2$.
For the coefficients of $x^{n-1}$ and $x^{n-2}$ in the general expansion
\eqref{polt} we have
\ben
\bal
\phi_1(u)&=E_{11}(u)+\dots+E_{n\ts n}(u),\\[0.2em]
\phi_2(u)&=\sum_{1\leqslant i<j\leqslant n}\ts\ts\begin{vmatrix}
E_{i\tss i}(u)&E_{i\tss j}(u)\\[0.4em]
E_{j\tss i}(u)&E_{j\tss j}(u)
                \end{vmatrix}+\sum_{i=1}^n\ts (i-1)\ts E^{\ts\prime}_{i\tss i}(u).
\eal
\een
Hence,
\ben
\bal
\phi^{(0)}_1&=E^{(0)}_{11}[-1]+\dots+E^{(0)}_{n\ts n}[-1]=e_{11}[-1]+\dots+e_{n+1\ts n+1}[-1],\\[0.3em]
\phi^{(1)}_1&=E^{(1)}_{n\tss n}[-1]=e_{n\ts n+1}[-1],
\eal
\een
whereas
\ben
\phi^{(1)}_2=\sum_{i=1}^{n-1}\ts\ts\begin{vmatrix}
E^{(0)}_{i\tss i}[-1]&E^{(1)}_{i\tss n}[-1]\\[0.4em]
E^{(0)}_{n\tss i}[-1]&E^{(1)}_{n\tss n}[-1]
                \end{vmatrix}+(n-1)\ts E^{(1)}_{n\tss n}[-2],
\een
cf. \cite[Sec.~5]{ap:qm}.
\qed
\ees

As with the minimal nilpotent case considered in \cite{ap:qm}, the proof
of Theorem~\ref{thm:suga} follows the
same approach as in the paper \cite{cm:ho} which deals with the case $e=0$.
We will give some necessary additional details
in Section~\ref{sec:pt} which include the use of the derivation $\De$
of the algebra $\U(\wh\agot)$ determined by its commutation relations with
the left multiplication operators
\beql{bde}
\big[\De,X[r]\tss\big]=r\ts X[r+1],\qquad
\big[\De,\vac\big]=0,
\eeq
for $X\in\agot$. We can regard $\De$ as an operator on the vacuum module $V(\agot)$.
It is immediate from the commutation relations that $\De$
preserves Segal--Sugawara vectors and thus gives rise to the operator
\ben
\De:\z(\wh\agot)\to \z(\wh\agot).
\een
Its use in the proof of Theorem~\ref{thm:suga} will be based on the following
property of the coefficients $\phi^{(r)}_k\in V(\agot)$, where we do not assume that they belong
to $\z(\wh\agot)$.

\ble\label{lem:acde}
Under the action of $\De$ we have
\beql{dephize}
\De:\phi^{(r)}_k\mapsto 0
\eeq
for all $r> \la_{n-k+2}+\dots+\la_n-k+1$, while for $r=\la_{n-k+2}+\dots+\la_n-k+1$
we have
\ben
\De:\phi^{(r)}_k\mapsto -(k-1)\tss(\la_1+\dots+\la_{n-k+1})\tss \phi^{(r)}_{k-1}.
\een
\ele

We will prove Lemma~\ref{lem:acde} in Section~\ref{sec:pt}. In particular, it implies that
a complete set of Segal--Sugawara vectors in the case $\la_1=\dots=\la_n=1$
(that is, for $\agot=\gl_n$) can be produced from a single vector. Set
\ben
\phi=\cdet\left[\begin{matrix}
T+E_{11}[-1]&E_{12}[-1]&\dots&E_{1n}[-1]\\[0.2em]
E_{21}[-1]&T+E_{22}[-1]&\dots&E_{2n}[-1]\\
\vdots&\vdots& \ddots&\vdots     \\
E_{n1}[-1]&E_{n2}[-1]&\dots&T+E_{n\tss n}[-1]
                \end{matrix}\right],
\een
which is the constant term of the polynomial in \eqref{cdete}; see also
\cite{cm:ho} and \cite[Ch.~7]{m:so}.

\bco\label{cor:ff}
The elements $\De^{k}\phi$ for $k=0,1,\dots,n-1$ form a complete set
of Segal--Sugawara vectors for $\gl_n$.
\qed
\eco

A slightly less explicit but equivalent expression for the coefficients $\phi^{(r)}_k$
is obtained in the following proposition, where we introduce $\tau$
as an element of the extended Lie algebra $\wh\agot\oplus\CC\tau$
which satisfies
the commutation relations
\ben
\big[\tau,X[r]\tss\big]=-r\ts X[r-1],\qquad
\big[\tau,\vac\big]=0.
\een
For all $i,j\in\{1,\dots,n\}$ define the elements $\Ec_{ij}(\tau)$
of the tensor product space $V(\agot)\ot\CC[\tau]$
by
\ben
\Ec_{ij}(\tau)=\begin{cases}\de_{ij}\tau^{\la_j}+E^{(0)}_{ij}[-1]\tss\tau^{\la_j-1}+\dots+E^{(\la_j-1)}_{ij}[-1]
\qquad&\text{if}\quad i\geqslant j,\\[0.4em]
E^{(\la_j-\la_i)}_{ij}[-1]\tss\tau^{\la_i-1}+\dots+E^{(\la_j-1)}_{ij}[-1]
\qquad&\text{if}\quad i< j.
\end{cases}
\een
Introduce elements $\phi^{\circ}_1,\dots,\phi^{\circ}_N\in V(\agot)$ by using the expansion of the
column-determinant of the matrix $\Ec(\tau)=[\Ec_{ij}(\tau)]$,
\ben
\cdet\ts \Ec(\tau)=\tau^N+\phi^{\circ}_1\ts\tau^{N-1}
+\dots+\phi^{\circ}_N.
\een
We will say that the element $E^{(b)}_{ij}[r]$ of
the Lie algebra $t^{-1}\agot[t^{-1}]$ has {\em weight $b$}.
This defines a grading
on the universal enveloping algebra
$\U\big(t^{-1}\agot[t^{-1}]\big)$.

\bpr\label{prop:sugaeq} For each $k=1,\dots,n$ and $r$ satisfying
\eqref{conda}, the Segal--Sugawara vector
$\phi^{(r)}_k$ coincides with the homogeneous component
of maximal weight of the coefficient $\phi^{\circ}_{r+k}$.
\epr

\bpf
Expand the
column-determinant
\ben
\cdet\ts \Ec(\tau)=\sum_{\si\in\Sym_n} \sgn\si\cdot \Ec_{\si(1)\tss 1}(\tau)\dots
\Ec_{\si(n)\tss n}(\tau)
\een
and move all powers of $\tau$ to the right by commuting them with the
elements $E^{(b)}_{ij}[r]$.
It is straightforward to verify that a nonzero contribution to the
maximum weight component of $\phi^{\circ}_{r+k}$ can only come from
commutators involving the leading terms $\tau^{\la_i}$ of $\Ec_{i\tss i}(\tau)$
for $i=1,\dots,n$. Moreover, this
component has weight equal to $r$ and coincides with $\phi^{(r)}_k$.
\epf

By adapting Rybnikov's construction \cite{r:si} to the case of
the affine Kac--Moody algebra $\wh\agot$ as in
\cite{ap:qm}, for any element $\chi\in\agot^*$ and a variable $z$
consider the homomorphism
\beql{evalrozmu}
\varrho^{}_{\tss\chi}:\U\big(t^{-1}\agot[t^{-1}]\big)\to \U(\agot)\ot\CC[z^{-1}],
\qquad X[r]\mapsto X\tss z^r+\de_{r,-1}\ts\chi(X),
\eeq
for any $X\in\agot$ and $r<0$.
If $S\in \z(\wh\agot)$ is a homogeneous
element of degree $d$ with respect to the grading defined by
$\deg X[r]=-r$ for $r<0$,
define the elements $S^{}_{(k)}\in\U(\agot)$
(depending on $\chi$) by the expansion
\ben
\varrho^{}_{\tss\chi}(S)=S^{}_{(0)}\ts z^{-d}+\dots+S^{}_{(d-1)}\ts z^{-1}+S^{}_{(d)}.
\een
If the variable $z$ takes
a particular nonzero value in $\CC$, then the formula in \eqref{evalrozmu}
defines a homomorphism $\U\big(t^{-1}\agot[t^{-1}]\big)\to \U(\agot)$.
Since $\z(\wh\agot)$ is a commutative
subalgebra of $\U\big(t^{-1}\agot[t^{-1}]\big)$, its image under
this homomorphism is a commutative subalgebra of $\U(\agot)$
which we denote by $\Ac_{\chi}$. This subalgebra
does not depend on the value of $z$.

It is clear from its definition that the degree of the coefficient $\phi^{(r)}_k$ in \eqref{polt}
equals $k$. Therefore, the respective degrees of the Segal--Sugawara vectors
provided by Theorem~\ref{thm:suga}
coincide with the
degrees of the basic invariants of the symmetric algebra $\Sr(\agot)$ given by
\ben
\underbrace{1,\dots,1}_{\la_n},\underbrace{2,\dots,2}_{\la_{n-1}},\dots,\underbrace{n,\dots,n}_{\la_1},
\een
as found in \cite{ppy:si}; see also \cite{bb:ei}.
Define elements $\phi^{(a)}_{k\tss (m)}\in\U(\agot)$ by
applying the homomorphism $\varrho^{}_{\tss\chi}$ to the Segal--Sugawara vectors:
\ben
\varrho^{}_{\tss\chi}: \phi^{(r)}_k\mapsto \phi^{(r)}_{k\tss (0)}\ts z^{-k}
+\dots+\phi^{(r)}_{k\tss (k-1)}\ts z^{-1}+\phi^{(r)}_{k\tss (k)}.
\een
More explicitly, set $\chi_{ij}^{(r)}=\chi(E_{ij}^{(r)})$ and
\ben
\chi_{ij}(u)=\begin{cases}\chi^{(0)}_{ij}+\dots+\chi^{(\la_j-1)}_{ij}\ts u^{\la_j-1}
&\text{if}\quad i\geqslant j,\\[0.4em]
\chi^{(\la_j-\la_i)}_{ij}\tss u^{\la_j-\la_i}+\dots+\chi^{(\la_j-1)}_{ij}\ts u^{\la_j-1}
&\text{if}\quad i< j.
\end{cases}
\een
Furthermore, let
\ben
\Ec_{ij}(u)=\begin{cases}E^{(0)}_{ij}+\dots+E^{(\la_j-1)}_{ij}\ts u^{\la_j-1}
&\text{if}\quad i\geqslant j,\\[0.4em]
E^{(\la_j-\la_i)}_{ij}\tss u^{\la_j-\la_i}+\dots+E^{(\la_j-1)}_{ij}\ts u^{\la_j-1}
&\text{if}\quad i< j.
\end{cases}
\een
The image of the column-determinant \eqref{cdete}
under the homomorphism $\varrho^{}_{\tss\chi}$ will be written as a polynomial in $x$,
\begin{multline}
\non
\cdet\left[\begin{matrix}
x-\la_1\tss\di_z+\chi_{11}(u)+\Ec_{11}(u)\tss z^{-1}&\dots&\chi_{1n}(u)+\Ec_{1n}(u)\tss z^{-1}\\
\vdots& \ddots&\vdots     \\
\chi_{n1}(u)+\Ec_{n1}(u)\tss z^{-1}&\dots&
x-\la_n\tss\di_z+\chi_{n\tss n}(u)+\Ec_{n\tss n}(u)\tss z^{-1}
                \end{matrix}\right]\\[1em]
{}=x^n+\phi_1(u,z)\tss x^{n-1}+\dots+\phi_n(u,z),
\end{multline}
with
\ben
\phi_k(u,z)=\sum_{r,m}\phi^{(r)}_{k\tss (m)}\ts u^r\tss z^{-k+m}.
\een

Recall that the {\em Mishchenko--Fomenko subalgebra} $\overline\Ac_{\chi}$
of the symmetric algebra $\Sr(\agot)$ is a Poisson-commutative subalgebra
generated by all $\chi$-shifts
of all $\agot$-invariants $P\in \Sr(\agot)^{\agot}$, as originally defined in
\cite{mf:ee}. By using the results of \cite{ap:qm} we come to the following.

\bco\label{cor:genamu}
Suppose the element $\chi\in\agot^*$ is regular.
Then the elements $\phi^{(r)}_{k\tss (m)}$
for $k=1,\dots,n$ with $r$ satisfying conditions \eqref{conda}
and $m=0,\dots,k-1$
are algebraically independent generators of the algebra $\Ac_{\chi}$. Moreover,
$\Ac_{\chi}$ is a quantization of the subalgebra
$\overline\Ac_{\chi}$
so that
$\gr\Ac_{\chi}=\overline\Ac_{\chi}$.
\qed
\eco

As another corollary of Theorem~\ref{thm:suga}, we recover the algebraically independent
generators of the center of $\U(\agot)$ constructed in \cite{bb:ei} with the use of the shifted
Yangians; see also \cite{m:ce}
for the particular case of rectangular pyramids. The algebra $\Ac_{\tss 0}$
coincides with the center of $\U(\agot)$ so that the formulas for the generators
are found by taking $\chi=0$ in \eqref{evalrozmu}; cf. \cite[Secs~6.5~\&~7.2]{m:so}.
Write the column-determinant
\ben
\cdet\left[\begin{matrix}
x+(n-1)\tss\la_1+\Ec_{11}(u)&\Ec_{12}(u)&\dots&\Ec_{1n}(u)\\[0.4em]
\Ec_{21}(u)&x+(n-2)\tss\la_2+\Ec_{22}(u)&\dots&\Ec_{2n}(u)\\
\vdots&\vdots& \ddots&\vdots     \\
\Ec_{n1}(u)&\Ec_{n2}(u)&\dots&x+\Ec_{n\tss n}(u)
                \end{matrix}\right]
\een
as a polynomial in $x$,
\ben
x^n+\Phi_1(u)\tss x^{n-1}+\dots+\Phi_n(u)\qquad\text{with}\quad
\Phi_k(u)=\sum_r\Phi^{(r)}_k\ts u^r.
\een

\bco\label{cor:gencent}
The coefficients $\Phi^{(r)}_k$ with $k=1,\dots,n$ and
$r$ satisfying conditions \eqref{conda}
are algebraically independent generators of the center of
the algebra $\U(\agot)$.
\qed
\eco

The elements $\Phi^{(r)}_k$ are slightly different from the corresponding central elements
$z_r$
given by \cite[Eq.~(1.3)]{bb:ei}. The mapping
$E_{ij}^{(p)}\mapsto E_{ij}^{(p)}+\de_{p\tss 0}\tss \de_{ij}\tss c\tss\la_i$
defines an automorphism of the algebra $\U(\agot)$ for
any given constant $c$. The image of the Casimir element $\Phi^{(r)}_k$
under this automorphism with $c=-n+1$ coincides with $z_{r+k}$.

\section{Proof of Theorem~\ref{thm:suga}}
\label{sec:pt}

To prove the first part of the theorem,
we will verify that the coefficients $\phi^{(r)}_k\in V(\agot)$ are annihilated
by a family of elements which generate $\agot[t]$ as a Lie algebra.
For such a family we take
$E^{(0)}_{i+1\ts i}[0]$, $E^{(\la_{i+1}-\la_i)}_{i\ts i+1}[0]$ for $i=1,\dots,n-1$,
and $E^{(p)}_{i\ts i}[s]$ for $p=0,\dots,\la_i-1$ all $i=1,\dots,n$ and $s\geqslant 0$.

We will apply the generators to
the column-determinant $\cdet\ts\Ec$ of the matrix $\Ec=[\Ec_{ij}]$, where we set
$\Ec_{ij}=\de_{ij}\tss(x+\la_i\tss T)+E_{ij}(u)$. The result of such an application
is
a polynomial in $x$ of the form
\beql{psipo}
\psi_1(u)\ts x^{n-1}
+\dots+\psi_n(u),\qquad \psi_k(u)=\sum_r\psi^{(r)}_k\ts u^r,\qquad
\psi^{(r)}_k\in V(\agot).
\eeq
It will be sufficient to verify that for $k=1,\dots,n$ the degree of the polynomial
$\psi_k(u)$ in $u$ is less than
\beql{lowb}
\la_{n-k+2}+\dots+\la_n-k+1
\eeq
(in particular, $\psi_1(u)$ must be equal to zero).
We will consider the groups of generators case by case and
use the notation $\tau_i=x+\la_i\tss T$ for brevity. Note
the commutation relations
\beql{tauri}
\big[\tau_i,X[r]\tss\big]=-r\ts\la_i\ts X[r-1],\qquad
X\in\agot.
\eeq

\paragraph{Generators $E^{(0)}_{i+1\ts i}[0]${\rm .}}

We will
rely on some simple properties of column-determinants described in
\cite[Lemmas~4.1 \& 4.2]{cm:ho}. They allow us to write
the element $E^{(0)}_{i+1\ts i}[0]\ts\cdet\ts\Ec$
as the difference of two column-determinants
\ben
\begin{vmatrix}
\Ec_{1\ts 1}&\dots&\Ec_{1\ts i}&\Ec_{1\ts i+1}&\dots&\Ec_{1\ts n}\\[0.2em]
\dots&\dots&\dots&\dots &\dots &\dots \\[0.2em]
\wt\Ec_{i+1\ts 1}&\dots&\wt\Ec_{i+1\ts i}&\wt\Ec_{i+1\ts i+1}&\dots&\wt\Ec_{i+1\ts n}\\[0.2em]
\Ec_{i+1\ts 1}&\dots&\Ec_{i+1\ts i}&\Ec_{i+1\ts i+1}&\dots&\Ec_{i+1\ts n}\\[0.2em]
\dots&\dots&\dots&\dots &\dots &\dots \\[0.2em]
\Ec_{n\ts 1}&\dots&\Ec_{n\ts i}&\Ec_{n\ts i+1}&\dots&\Ec_{n\ts n}
                \end{vmatrix}
\quad - \quad
                \begin{vmatrix}
\Ec_{1\ts 1}&\dots&\Ec_{1\ts i}&\wh\Ec_{1\ts i}&\dots&\Ec_{1\ts n}\\[0.2em]
\dots&\dots&\dots&\dots &\dots &\dots \\[0.2em]
\Ec_{i\ts 1}&\dots&\Ec_{i\ts i}&\wh\Ec_{i\ts i}&\dots&\Ec_{i\ts n}\\[0.2em]
\Ec_{i+1\ts 1}&\dots&\Ec_{i+1\ts i}&\wh\Ec_{i+1\ts i}&\dots&\Ec_{i+1\ts n}\\[0.2em]
\dots&\dots&\dots&\dots &\dots &\dots \\[0.2em]
\Ec_{n\ts 1}&\dots&\Ec_{n\ts i}&\wh\Ec_{n\ts i}&\dots&\Ec_{n\ts n}
                \end{vmatrix},
\een
obtained from $\cdet\ts\Ec$ by replacing row $i$ and column $i+1$, respectively, as indicated.
Here we set $\wt\Ec_{i+1\ts j}=\Ec_{i+1\ts j}$ for $j\leqslant i$ and
\ben
\wt\Ec_{i+1\ts j}=\de^{}_{i+1\ts j}\tau^{}_{i+1}+E^{(\la_j-\la_i)}_{i+1\ts j}[-1]\tss u^{\la_j-\la_i}
+\dots+E^{(\la_j-1)}_{i+1\ts j}[-1]\tss u^{\la_j-1}
\een
for $j\geqslant i+1$, while $\wh\Ec_{k\tss i}=\Ec_{k\tss i}$ for $k\geqslant i+1$ and
\ben
\wh\Ec_{k\tss i}=
\de^{}_{k\tss i}\tss\tau^{}_{i+1}+E^{(\la_{i+1}-\la_k)}_{k\tss i}[-1]\tss u^{\la_{i+1}-\la_k}
+\dots+E^{(\la_i-1)}_{k\tss i}[-1]\tss u^{\la_i-1}
\een
for $k\leqslant i$.

Start with
the first determinant and observe that it stays unchanged if
we subtract row $i+1$ from row $i$. The first $i$ entries in the new row $i$ will be
equal to zero and it will have the form
\ben
[0,\dots,0,\check{\Ec}_{i+1\ts i+1},\dots,\check{\Ec}_{i+1\ts n}]
\een
with
\ben
\check{\Ec}_{i+1\ts j}=
E^{(\la_j-\la_{i+1})}_{i+1\ts j}[-1]\tss u^{\la_j-\la_{i+1}}+\dots+
E^{(\la_j-\la_i-1)}_{i+1\ts j}[-1]\tss u^{\la_j-\la_i-1}
\een
for $j=i+1,\dots,n$. Write the first determinant as a polynomial in $x$
and consider the coefficient of $x^{n-k}$. This coefficient is a polynomial in $u$,
and if $n-k\leqslant i-1$, then
its degree does not exceed
\beql{sumi}
(\la_{n-k+1}-1)+\dots+(\la_j-\la_i-1)+\dots+(\la_n-1).
\eeq
This is clear since $E_{i\tss l}(u)$ is a polynomial
in $u$ of degree $\la_l-1$, while $\check{\Ec}_{i+1\ts j}$ is
a polynomial of degree $\la_j-\la_i-1$. However, the sum \eqref{sumi}
is less than \eqref{lowb} because $\la_{n-k+1}-\la_i-1<0$.

Similarly, if $n-k\geqslant i$, then the coefficient of
$x^{n-k}$ is a polynomial in $u$ whose degree does not exceed the sum
\ben
(\la_{n-k+2}-1)+\dots+(\la_j-\la_i-1)+\dots+(\la_n-1)
\een
which is less than \eqref{lowb} because $\la_i>0$. Therefore, the contribution
to the polynomial $\psi_k(u)$ in \eqref{psipo} arising from the first determinant
does not violate the required degree condition.

To reach the same conclusion for the second determinant,
use its simultaneous expansion along columns $i$ and $i+1$.
The expansion will involve
$2\times 2$ column-minors of the form
\beql{minortt}
\begin{vmatrix}
\Ec_{a\tss i}&\wh{\Ec}_{a\tss i}\\[0.4em]
\Ec_{b\tss i}&\wh{\Ec}_{b\tss i}
\end{vmatrix}
\eeq
for $a<b$. First consider the case with $a=i$
and suppose that $\la_i=\la_{i+1}$. Then the minor \eqref{minortt} equals
\ben
\big(\tau_i+E_{i\tss i}(u)\big)E_{b\tss i}(u)-E_{b\tss i}(u)\big(\tau_{i+1}+E_{i\tss i}(u)\big)
=E_{b\tss i}(u)(\tau_i-\tau_{i+1})+\big[\tau_i+E_{i\tss i}(u),E_{b\tss i}(u)\big].
\een
We have $E_{b\tss i}(u)(\tau_i-\tau_{i+1})=E_{b\tss i}(u)(\la_i-\la_{i+1})\tss T=0$.
Note that since
$[E^{(r)}_{i\tss i},E^{(s)}_{b\tss i}]=0$ for the values
with the condition $r+s\geqslant \la_i$, the degree of the polynomial
$\big[\tau_i+E_{i\tss i}(u),E_{b\tss i}(u)\big]$ in $u$ does not exceed $\la_i-1$.
Furthermore, by \eqref{comm} and \eqref{tauri} the coefficient of $u^{\la_i-1}$ equals
$\la_i E^{(\la_i-1)}_{b\tss i}[-2]-\la_i E^{(\la_i-1)}_{b\tss i}[-2]=0$ so that the degree
of the polynomial does not exceed $\la_i-2$.

Now consider the contribution of the $2\times 2$ minor under consideration to
the coefficient of $x^{n-k}$ in the second determinant.
If $n-k\leqslant i-1$ then the contribution
is a polynomial in $u$ whose degree
does not exceed
\ben
(\la_{n-k+1}-1)+\dots+(\la_{i-1}-1)+(\la_{i}-2)+(\la_{i+2}-1)+\dots+(\la_n-1)
\een
which is less than the sum in \eqref{lowb} because $\la_{n-k+1}-1<\la_{i+1}$.
If $n-k\geqslant i$ then the contribution to
the coefficient of $x^{n-k}$ is a polynomial in $u$ whose degree
does not exceed
\ben
(\la_{i}-2)+(\la_{n-k+3}-1)+\dots+(\la_n-1)
\een
which is again less than \eqref{lowb}.

Continuing with the case $a=i$ in \eqref{minortt}, suppose now that $\la_i<\la_{i+1}$.
Exactly as above we find that the minor equals $E_{b\tss i}(u)(\tau_i-\tau_{i+1})$
plus a polynomial in $u$ of degree $\leqslant\la_{i+1}-2$.
This implies that the resulting contribution to
the coefficient of $x^{n-k}$ is a polynomial in $u$ whose degree
is less than \eqref{lowb}.

The minor \eqref{minortt} is equal to zero for $i+1\leqslant a<b$, while
the remaining values of $a$ and $b$ are considered in a way quite similar
to the above arguments. This allows us to conclude that the coefficients $\phi^{(r)}_k\in V(\agot)$
satisfying \eqref{conda}
are annihilated
by all generators of the form
$E^{(0)}_{i+1\ts i}[0]$.

\paragraph{Generators $E^{(\la_{i+1}-\la_i)}_{i\ts i+1}[0]${\rm .}}
We will argue in the same way as above.
It will be convenient to multiply the expression
$E^{(\la_{i+1}-\la_i)}_{i\ts i+1}[0]\ts\cdet\ts\Ec$ by $u^{\la_{i+1}-\la_i}$
so that the product will be written
as the difference of two column-determinants
\ben
\begin{vmatrix}
\Ec_{1\ts 1}&\dots&\Ec_{1\ts i}&\Ec_{1\ts i+1}&\dots&\Ec_{1\ts n}\\[0.2em]
\dots&\dots&\dots&\dots &\dots &\dots \\[0.2em]
\Ec_{i\ts 1}&\dots&\Ec_{i\tss i}&\Ec_{i\ts i+1}&\dots&\Ec_{i\tss  n}\\[0.2em]
\wt\Ec_{i\tss 1}&\dots&\wt\Ec_{i\tss i}&\wt\Ec_{i\ts i+1}&\dots&\wt\Ec_{i\tss n}\\[0.2em]
\dots&\dots&\dots&\dots &\dots &\dots \\[0.2em]
\Ec_{n\ts 1}&\dots&\Ec_{n\ts i}&\Ec_{n\ts i+1}&\dots&\Ec_{n\ts n}
                \end{vmatrix}
\quad - \quad
                \begin{vmatrix}
\Ec_{1\ts 1}&\dots&\wh\Ec_{1\ts i+1}&\Ec_{1\ts i+1}&\dots&\Ec_{1\ts n}\\[0.2em]
\dots&\dots&\dots&\dots &\dots &\dots \\[0.2em]
\Ec_{i\ts 1}&\dots&\wh\Ec_{i\ts i+1}&\Ec_{i\ts i+1}&\dots&\Ec_{i\ts n}\\[0.2em]
\Ec_{i+1\ts 1}&\dots&\wh\Ec_{i+1\ts i+1}&\Ec_{i+1\ts i+1}&\dots&\Ec_{i+1\ts n}\\[0.2em]
\dots&\dots&\dots&\dots &\dots &\dots \\[0.2em]
\Ec_{n\ts 1}&\dots&\wh\Ec_{n\ts i+1}&\Ec_{n\ts i+1}&\dots&\Ec_{n\ts n}
                \end{vmatrix},
\een
obtained from $\cdet\ts\Ec$ by replacing row $i+1$ and column $i$, respectively, as indicated,
where now we use the notation
$\wt\Ec_{i\tss j}=\Ec_{i\tss j}$ for $j\geqslant i+1$ and
\ben
\wt\Ec_{i\tss j}=\de^{}_{i\tss j}\tau^{}_{i}+E^{(\la_{i+1}-\la_i)}_{i\tss j}[-1]\tss u^{\la_{i+1}-\la_i}
+\dots+E^{(\la_j-1)}_{i\tss j}[-1]\tss u^{\la_j-1}
\een
for $j\leqslant i$, while $\wh\Ec_{k\ts i+1}=\Ec_{k\ts i+1}$ for $k\leqslant i$ and
\ben
\wh\Ec_{k\ts i+1}=
\de^{}_{k\ts i+1}\tss\tau^{}_{i}+E^{(\la_{i+1}-\la_i)}_{k\ts i+1}[-1]\tss u^{\la_{i+1}-\la_i}
+\dots+E^{(\la_{i+1}-1)}_{k\ts i+1}[-1]\tss u^{\la_{i+1}-1}
\een
for $k\geqslant i+1$. The product
$u^{\la_{i+1}-\la_i}\ts E^{(\la_{i+1}-\la_i)}_{i\ts i+1}[0]\ts\cdet\ts\Ec$
is a polynomial in $x$ and we use calculations similar to the
previous case to verify that
for $k=1,\dots,n$ its coefficient of $x^{n-k}$ is a polynomial in $u$
of degree less than the sum of $\la_{i+1}-\la_i$ and the expression \eqref{lowb}.

\paragraph{Generators $E^{(0)}_{i\tss i}[1]${\rm .}}
Note the commutation relations $\big[E^{(0)}_{i\tss i}[1],\tau_j\big]=\la_j\tss E^{(0)}_{i\tss i}[0]$
implied by \eqref{tauri}. Hence, we also have
the relations
\ben
\big[E^{(0)}_{i\tss i}[1], \Ec_{m\tss l}\big]=
\de_{m\tss l}\tss\big(\la_m\tss E^{(0)}_{i\tss i}[0]+\vk_{i\tss m}\big)
+\de_{m\tss i}\tss \Ec_{i\tss l}[0]-\de_{i\tss l}\tss \Ec_{m\tss i}[0],
\een
where we set $\vk_{i\tss m}=\big\langle E_{ii}^{(0)},E_{mm}^{(0)}\big\rangle$
and use the notation
\beql{eijze}
\Ec_{ij}[0]=\begin{cases}E^{(0)}_{ij}[0]+\dots+E^{(\la_j-1)}_{ij}[0]\ts u^{\la_j-1}
&\text{if}\quad i\geqslant j,\\[0.4em]
E^{(\la_j-\la_i)}_{ij}[0]\tss u^{\la_j-\la_i}+\dots+E^{(\la_j-1)}_{ij}[0]\ts u^{\la_j-1}
&\text{if}\quad i< j.
\end{cases}
\eeq
By \cite[Lemma~4.1]{cm:ho}, we get
\beql{eoneact}
\big[E^{(0)}_{i\tss i}[1],\cdet\ts\Ec\big]=
\sum_{j=1}^n\ts
\begin{vmatrix}
\Ec_{1\tss 1}&\dots&\big[E^{(0)}_{i\tss i}[1], \Ec_{1\tss j}\big]&\dots&\Ec_{1\tss n}\\[0.2em]
\dots&\dots&\dots&\dots &\dots \\[0.2em]
\Ec_{n\tss 1}&\dots&\big[E^{(0)}_{i\tss i}[1], \Ec_{n\tss j}\big]&\dots&\Ec_{n\tss n}
                \end{vmatrix}
\eeq
which equals the sum
\beql{sumet}
\sum_{j=1}^n\tss\vk_{ij}\ts \Ec^{\wh\jmath}_{\wh\jmath}+
\sum_{j=1}^n\ts
\begin{vmatrix}
\Ec_{1\tss 1}&\dots&0&\dots&\Ec_{1\tss n}\\[0.2em]
\dots&\dots&\dots&\dots &\dots \\[0.2em]
\dots&\dots&\la_j\ts E^{(0)}_{i\tss i}[0] &\dots &\dots \\[0.2em]
\dots&\dots&\dots&\dots &\dots \\[0.2em]
\Ec_{n\tss 1}&\dots&0&\dots&\Ec_{n\tss n}
                \end{vmatrix}
\eeq
plus the difference of two column-determinants
\beql{difftwo}
\begin{vmatrix}
\Ec_{1\tss 1}&\dots&\Ec_{1\tss n}\\[0.2em]
\dots&\dots&\dots \\[0.2em]
\Ec_{i\tss 1}[0]&\dots&\Ec_{i\tss n}[0] \\[0.2em]
\dots&\dots&\dots \\[0.2em]
\Ec_{n\tss 1}&\dots&\Ec_{n\tss n}
                \end{vmatrix}
\quad - \quad
                \begin{vmatrix}
\Ec_{1\tss 1}&\dots&\Ec_{1\ts i}[0]&\dots&\Ec_{1\tss n}\\[0.2em]
\dots&\dots&\dots&\dots &\dots \\[0.2em]
\Ec_{n\tss 1}&\dots&\Ec_{n\ts i}[0]&\dots&\Ec_{n\tss n}
                \end{vmatrix},
\eeq
where $\Ec^{\wh\jmath}_{\wh\jmath}$ denotes the column-determinant of the matrix
obtained from $\Ec$ by deleting row and column $j$.
Now we proceed as in \cite{cm:ho} relying on Lemma~4.2 therein to evaluate the action of
the elements of the form $E^{(0)}_{i\tss i}[0]$, $\Ec_{i\tss j}[0]$ and $\Ec_{m\tss i}[0]$.
For the generator $E^{(0)}_{i\tss i}[0]$ occurring as the $(j,j)$-entry with $j\ne i$
we have
\beql{midt}
\begin{vmatrix}
\Ec_{1\tss 1}&\dots&0&\dots&\Ec_{1\tss n}\\[0.2em]
\dots&\dots&\dots&\dots &\dots \\[0.2em]
\dots&\dots&E^{(0)}_{i\tss i}[0] &\dots &\dots \\[0.2em]
\dots&\dots&\dots&\dots &\dots \\[0.2em]
\Ec_{n\tss 1}&\dots&0&\dots&\Ec_{n\tss n}
\end{vmatrix}
=
\sum_{m=j+1}^n\ts
\begin{vmatrix}
\Ec_{1\tss 1}&\dots&0&\dots&\Ec_{1\tss n}\\[0.2em]
\dots&\dots&\dots&\dots &\dots \\[0.2em]
\dots&\dots&\Ec_{i\tss m} &\dots &\dots \\[0.2em]
\dots&\dots&\dots&\dots &\dots \\[0.2em]
\Ec_{n\tss 1}&\dots&0&\dots&\Ec_{n\tss n}
\end{vmatrix}-\de^{}_{j<i}\ts\Ec^{\wh\jmath}_{\wh\jmath}\ts,
\eeq
where row and column $j$ are deleted in the column-determinants in the sums, and
$\Ec_{i\tss m}$ occurs in row $i$ and column $m$. Clearly,
the left hand side of \eqref{midt}
is zero for $j=i$.

Next consider the first determinant in \eqref{difftwo}.
For the expression $\Ec_{i\tss j}[0]$ occurring as the $(i,j)$ entry,
with the remaining entries in column $j$ equal to zero,
we get an expansion into a sum of column-determinants,
\beql{eact}
\begin{vmatrix}
\Ec_{1\tss 1}&\dots&0&\dots&\Ec_{1\tss n}\\[0.2em]
\dots&\dots&\dots&\dots &\dots \\[0.2em]
\dots&\dots&\Ec_{i\tss j}[0] &\dots &\dots \\[0.2em]
\dots&\dots&\dots&\dots &\dots \\[0.2em]
\Ec_{n\tss 1}&\dots&0&\dots&\Ec_{n\tss n}
\end{vmatrix}
=
(-1)^{i+j}\sum_{m=j+1}^n\ts
\begin{vmatrix}
\Ec_{1\tss 1}&\dots&\big[\Ec_{i\tss j}[0],\Ec_{1\tss m}\big]&\dots&\Ec_{1\tss n}\\[0.2em]
\dots&\dots&\dots&\dots &\dots \\[0.2em]
\Ec_{n\tss 1}&\dots&\big[\Ec_{i\tss j}[0],\Ec_{n\tss m}\big]&\dots&\Ec_{n\tss n}
\end{vmatrix},
\eeq
where row $i$ and column $j$ are deleted in the column-determinants in the sum.
Note that nonzero commutators occur only for $i\ne j$. If $m\ne i$, then the
commutator $\big[\Ec_{i\tss j}[0],\Ec_{j\tss m}\big]$ is a polynomial in $u$ whose
component of maximal degree
equals $\min(\la_i,\la_j)\tss E^{(\la_m-1)}_{i\tss m}[-1]\ts u^{\la_m-1}$.
By writing the column-determinants
as polynomials in $x$, we can see that the resulting contribution
to each polynomial $\psi_k(u)$ in
\eqref{psipo} of any
components with the powers of $u$
less than $\la_m-1$, is a polynomial in $u$ whose degree is less than \eqref{lowb}.
For this reason we may replace
the commutator $\big[\Ec_{i\tss j}[0],\Ec_{j\tss m}\big]$ in the above expansion
with the polynomial $\min(\la_i,\la_j)\tss \Ec_{i\tss m}$.

Furthermore, observe that if $s\ne j$ and $s<i$, then the
commutator $\big[\Ec_{i\tss j}[0],\Ec_{s\tss i}\big]$
occurring for $i>j$,
is zero for $\la_i-\la_s\geqslant\la_j$.
Otherwise, if $\la_i-\la_s<\la_j$, then the
commutator is a polynomial in $u$ whose
component of maximal degree
equals $(\la_i-\la_j-\la_s)\tss E^{(\la_j-1)}_{s\tss j}[-1]\ts u^{\la_j-1}$.
Similarly, the maximal degree component of the commutator
$\big[\Ec_{i\tss j}[0],\Ec_{s\tss i}\big]$ with $s>i>j$ equals
$-\la_j\tss E^{(\la_j-1)}_{s\tss j}[-1]\ts u^{\la_j-1}$.

By looking at the powers of $x$ we find that
the contributions of these components to each polynomial $\psi_k(u)$ in
\eqref{psipo} are polynomials in $u$ whose degrees are less than the expression in \eqref{lowb},
unless $\la_j=\la_i$.
Moreover, the same property is shared by the components with
powers of $u$ less than $\la_j-1$. This allows us
to replace all commutators $\big[\Ec_{i\tss j}[0],\Ec_{s\tss i}\big]$ with $s\ne j$
occurring in the column-determinant,
with the polynomial $-\min(\la_j,\la_s)\tss \Ec_{s\tss j}$.

The same argument implies that we may replace the remaining commutators with $i>j$ by the rule
\ben
\big[\Ec_{i\tss j}[0],\Ec_{j\tss i}\big]\rightsquigarrow \la_j(\Ec_{i\tss i}-\tau_i)
-\la_j(\Ec_{j\tss j}-\tau_j).
\een

Bringing the calculations together, we can write the first column-determinant in \eqref{difftwo}
as the sum of two expressions
\beql{sumon}
\sum_{j=1,\ j\ne i}^{n-1}(-1)^{i+j}\min(\la_i,\la_j)\sum_{m=j+1}^n
\begin{vmatrix}
\Ec_{1\tss 1}&\dots&0&\dots&\Ec_{1\tss n}\\[0.2em]
\dots&\dots&\dots&\dots &\dots \\[0.2em]
\dots&\dots&\Ec_{i\tss m}-\de_{i\tss m}\tau_i &\dots &\dots \\[0.2em]
\dots&\dots&\dots&\dots &\dots \\[0.2em]
\Ec_{n\tss 1}&\dots&0&\dots&\Ec_{n\tss n}
\end{vmatrix}
\eeq
and
\beql{sumtw}
-\sum_{j=1}^{i-1}(-1)^{i+j}
\begin{vmatrix}
\Ec_{1\tss 1}&\dots&\la_1\tss \Ec_{1\tss j}&\dots&\Ec_{1\tss n}\\[0.2em]
\dots&\dots&\dots&\dots &\dots \\[0.2em]
\dots&\dots&\la_j\tss (\Ec_{j\tss j}-\tau_j) &\dots &\dots \\[0.2em]
\dots&\dots&\dots&\dots &\dots \\[0.2em]
\Ec_{n\tss 1}&\dots&\la_j\tss \Ec_{n\tss j}&\dots&\Ec_{n\tss n}
\end{vmatrix}.
\eeq
Here row $i$ and column $j$ are considered to be deleted in the
column-determinants so that the rows and columns are labelled,
respectively, by the symbols $1,\dots,\wh\imath,\dots,n$ and
$1,\dots,\wh\jmath,\dots,n$. The only nonzero entry
$\Ec_{i\tss m}-\de_{i\tss m}\tau_i$ in column $m$ in the
column-determinants in \eqref{sumon} occurs in row $j$, while
the displayed middle column in \eqref{sumtw} occupies column $i$.

Note that permutation of rows in column-determinants results in a changed
sign as for usual commutative determinants. Hence, by moving row $j$
in the expression \eqref{sumon} to the position of the $i$-th row we can
present this expression in the form
\beql{sumonmo}
-\sum_{j=1,\ j\ne i}^{n-1}\min(\la_i,\la_j)\sum_{m=j+1}^n
\begin{vmatrix}
\Ec_{1\tss 1}&\dots&0&\dots&\Ec_{1\tss n}\\[0.2em]
\dots&\dots&\dots&\dots &\dots \\[0.2em]
\dots&\dots&\Ec_{i\tss m}-\de_{i\tss m}\tau_i &\dots &\dots \\[0.2em]
\dots&\dots&\dots&\dots &\dots \\[0.2em]
\Ec_{n\tss 1}&\dots&0&\dots&\Ec_{n\tss n}
\end{vmatrix},
\eeq
where the entry $\Ec_{i\tss m}-\de_{i\tss m}\tau_i$ now occupies the $(i,m)$-position,
while the rows and columns of the column-determinant
are labelled by the symbols $1,\dots,\wh\jmath,\dots,n$.

In the expression \eqref{sumtw} write each column-determinant as the difference
\beql{sumtwdi}
\begin{vmatrix}
\Ec_{1\tss 1}&\dots&\la_1\tss \Ec_{1\tss j}&\dots&\Ec_{1\tss n}\\[0.2em]
\dots&\dots&\dots&\dots &\dots \\[0.2em]
\dots&\dots&\la_j\tss \Ec_{j\tss j} &\dots &\dots \\[0.2em]
\dots&\dots&\dots&\dots &\dots \\[0.2em]
\Ec_{n\tss 1}&\dots&\la_j\tss \Ec_{n\tss j}&\dots&\Ec_{n\tss n}
\end{vmatrix}
\quad-\quad
\begin{vmatrix}
\Ec_{1\tss 1}&\dots&0&\dots&\Ec_{1\tss n}\\[0.2em]
\dots&\dots&\dots&\dots &\dots \\[0.2em]
\dots&\dots&\la_j\tss \tau_j &\dots &\dots \\[0.2em]
\dots&\dots&\dots&\dots &\dots \\[0.2em]
\Ec_{n\tss 1}&\dots&0&\dots&\Ec_{n\tss n}
\end{vmatrix}
\eeq
and move row $j$ in the second determinant to the position of the $i$-th row
so that the corresponding part of \eqref{sumtw} will be written as the sum
\ben
-\sum_{j=1}^{i-1}\
\begin{vmatrix}
\Ec_{1\tss 1}&\dots&0&\dots&\Ec_{1\tss n}\\[0.2em]
\dots&\dots&\dots&\dots &\dots \\[0.2em]
\dots&\dots&\la_j\tss \tau_j &\dots &\dots \\[0.2em]
\dots&\dots&\dots&\dots &\dots \\[0.2em]
\Ec_{n\tss 1}&\dots&0&\dots&\Ec_{n\tss n}
\end{vmatrix},
\een
where the entry $\la_j\tss \tau_j$ now occupies the $(i,i)$-position.

The next step is to show that taking signs into account, we can move column $i$
in the first column-determinant
in \eqref{sumtwdi} to the left by permuting it with
columns $i-1,i-2,\dots,j+1$ consecutively without changing its contribution
to the coefficients of the powers of $x$. Arguing by induction, suppose that column $i$
was swapped with a few columns on its left. As an induction step, we need
to show that for any $l\in\{j+1,\dots,i-1\}$ the sum of two column determinants of the form
\beql{difin}
\begin{vmatrix}
\Ec_{1\tss 1}&\dots&\Ec_{1\tss l}&\la_1\tss \Ec_{1\tss j}&\dots&\Ec_{1\tss n}\\[0.2em]
\dots&\dots&\dots&\dots&\dots &\dots \\[0.2em]
\dots&\dots&\Ec_{j\tss l} &\la_j\tss \Ec_{j\tss j} &\dots &\dots \\[0.2em]
\dots&\dots&\dots&\dots&\dots &\dots \\[0.2em]
\Ec_{n\tss 1}&\dots&\Ec_{n\tss l}&\la_j\tss \Ec_{n\tss j}&\dots&\Ec_{n\tss n}
\end{vmatrix}
\quad+\quad
\begin{vmatrix}
\Ec_{1\tss 1}&\dots&\la_1\tss \Ec_{1\tss j}&\Ec_{1\tss l}&\dots&\Ec_{1\tss n}\\[0.2em]
\dots&\dots&\dots&\dots&\dots &\dots \\[0.2em]
\dots&\dots&\la_j\tss \Ec_{j\tss j}&\Ec_{j\tss l} &\dots &\dots \\[0.2em]
\dots&\dots&\dots&\dots&\dots &\dots \\[0.2em]
\Ec_{n\tss 1}&\dots&\la_j\tss \Ec_{n\tss j}&\Ec_{n\tss l}&\dots&\Ec_{n\tss n}
\end{vmatrix}
\eeq
does not contribute to any polynomial $\psi_k(u)$ in
\eqref{psipo} beyond powers of $u$ less than \eqref{lowb}.
Using their expansions along the two displayed
adjacent columns we note that the sum of any $2\times 2$ minors
\beql{tbt}
\begin{vmatrix}
\Ec_{a\tss l}&\min(\la_a,\la_j)\tss \Ec_{a\tss j}\\
\Ec_{b\tss l}&\min(\la_b,\la_j)\tss \Ec_{b\tss j}
\end{vmatrix}
\quad+\quad
\begin{vmatrix}
\min(\la_a,\la_j)\tss \Ec_{a\tss j}&\Ec_{a\tss l}\\
\min(\la_b,\la_j)\tss \Ec_{b\tss j}&\Ec_{b\tss l}
\end{vmatrix}
\eeq
with $a<b$ is a polynomial in $u$ whose degree does not exceed $\la_l-1$.
Therefore, the contribution of the sum \eqref{difin} to $\psi_k(u)$
is a polynomial in $u$ whose degree is less than
the sum in \eqref{lowb}, except for the case where $\la_i=1$.
However, in this case we must have $\la_1=\dots=\la_i$ so that
the difference \eqref{tbt} is easily checked to be zero
(which agrees with the Manin matrix property of \cite[Lemma~4.3]{cm:ho}).
Thus, expression \eqref{sumtw} can be replaced with
\beql{finon}
\sum_{j=1}^{i-1}\
\begin{vmatrix}
\Ec_{1\tss 1}&\dots&\la_1\tss \Ec_{1\tss j}&\dots&\Ec_{1\tss n}\\[0.2em]
\dots&\dots&\dots&\dots &\dots \\[0.2em]
\dots&\dots&\la_j\tss \Ec_{j\tss j} &\dots &\dots \\[0.2em]
\dots&\dots&\dots&\dots &\dots \\[0.2em]
\Ec_{n\tss 1}&\dots&\la_j\tss \Ec_{n\tss j}&\dots&\Ec_{n\tss n}
\end{vmatrix}
\quad-\quad
\sum_{j=1}^{i-1}\
\begin{vmatrix}
\Ec_{1\tss 1}&\dots&0&\dots&\Ec_{1\tss n}\\[0.2em]
\dots&\dots&\dots&\dots &\dots \\[0.2em]
\dots&\dots&\la_j\tss \tau_j &\dots &\dots \\[0.2em]
\dots&\dots&\dots&\dots &\dots \\[0.2em]
\Ec_{n\tss 1}&\dots&0&\dots&\Ec_{n\tss n}
\end{vmatrix},
\eeq
where the middle column in the first sum
occupies the position of column $j$, with
the row and column $i$
considered to be deleted, while the middle column
occurring in the second sum
occupies the position of column $i$ with
the row and column $j$ considered to be deleted.

Observe that the column expansion along column $i$ and the use of
relations of the form \eqref{eact} lead to a similar evaluation of
the second determinant in \eqref{difftwo}. As in the above arguments we will only keep
the terms which can contribute to the coefficients of the powers $x^{n-k}$
beyond powers of $u$ less than \eqref{lowb}.
By an appropriate adjustment of the labels of rows and columns of emerging column-determinants,
as compared with the evaluation of the first determinant in \eqref{difftwo},
we can write the second determinant (taking into account
the minus sign in front) as the sum of two expressions
\beql{sumtmo}
\sum_{j=i+1}^n\
\begin{vmatrix}
\Ec_{1\tss 1}&\dots&\la_1\tss \Ec_{1\tss j}&\dots&\Ec_{1\tss n}\\[0.2em]
\dots&\dots&\dots&\dots &\dots \\[0.2em]
\dots&\dots&\la_{i-1}\tss \Ec_{i-1\ts j} &\dots &\dots \\[0.2em]
\dots&\dots&\la_i\tss \Ec_{i+1\ts j} &\dots &\dots \\[0.2em]
\dots&\dots&\dots&\dots &\dots \\[0.2em]
\dots&\dots&\la_i\tss (\Ec_{j\ts j}-\tau_j) &\dots &\dots \\[0.2em]
\dots&\dots&\dots&\dots &\dots \\[0.2em]
\Ec_{n\tss 1}&\dots&\la_i\tss \Ec_{n\tss j}&\dots&\Ec_{n\tss n}
\end{vmatrix}
\eeq
and
\beql{decse}
-\sum_{j=i+1}^{n}\
\begin{vmatrix}
\Ec_{1\tss 1}&\dots&\la_1\tss \Ec_{1\tss i}&\dots&\Ec_{1\tss n}\\[0.2em]
\dots&\dots&\dots&\dots &\dots \\[0.2em]
\dots&\dots&\la_i\tss \Ec_{i\tss i} &\dots &\dots \\[0.2em]
\dots&\dots&\dots&\dots &\dots \\[0.2em]
\Ec_{n\tss 1}&\dots&\la_i\tss \Ec_{n\tss i}&\dots&\Ec_{n\tss n}
\end{vmatrix}
\quad+\quad
\sum_{j=i+1}^{n}\
\begin{vmatrix}
\Ec_{1\tss 1}&\dots&0&\dots&\Ec_{1\tss n}\\[0.2em]
\dots&\dots&\dots&\dots &\dots \\[0.2em]
\dots&\dots&\la_i\tss \tau_i &\dots &\dots \\[0.2em]
\dots&\dots&\dots&\dots &\dots \\[0.2em]
\Ec_{n\tss 1}&\dots&0&\dots&\Ec_{n\tss n}
\end{vmatrix},
\eeq
where the rows and columns of
the column-determinants in \eqref{sumtmo} and in the second sum in \eqref{decse}
are labelled by $1,\dots,\wh\imath,\dots,n$ with $\la_i\tss \tau_i$ occurring as the $(j,j)$-entry,
while the rows and columns in the first sum in \eqref{decse}
are labelled by $1,\dots,\wh\jmath,\dots,n$.

Bringing together the evaluations of the column-determinants
in \eqref{sumet} and \eqref{difftwo} by taking into account \eqref{midt},
\eqref{sumonmo}, \eqref{finon}, \eqref{sumtmo} and \eqref{decse},
we can conclude that, modulo terms not contributing to the coefficients
of the powers $x^{n-k}$ beyond powers of $u$ less than \eqref{lowb},
the element $E^{(0)}_{i\tss i}[1]\tss\cdet\ts\Ec$
equals the sum of the following expressions:
\beql{explam}
\sum_{j=1}^{i-1}\
\begin{vmatrix}
\ \ \dots&(\la_1-\la_j)\tss \Ec_{1\tss j}&\dots\ \ \\[0.2em]
\ \ \dots&\dots&\dots\ \  \\[0.2em]
\ \ \dots&(\la_{j-1}-\la_j)\tss \Ec_{j-1\ts j} &\dots\ \  \\[0.2em]
\ \ \dots&0&\dots\ \ \\[0.2em]
\ \ \dots&\dots&\dots\ \  \\[0.2em]
\ \ \dots&0&\dots\ \
\end{vmatrix}
\quad+\quad
\sum_{j=i+1}^{n}\
\begin{vmatrix}
\ \ \dots&(\la_1-\la_i)\tss \Ec_{1\tss j}&\dots\ \ \\[0.2em]
\ \ \dots&\dots&\dots\ \  \\[0.2em]
\ \ \dots&(\la_{i-1}-\la_i)\tss \Ec_{j-1\ts j} &\dots\ \  \\[0.2em]
\ \ \dots&0&\dots\ \ \\[0.2em]
\ \ \dots&\dots&\dots\ \  \\[0.2em]
\ \ \dots&0&\dots\ \
\end{vmatrix},
\eeq
where the rows and columns of
the column-determinants are labelled by $1,\dots,\wh\imath,\dots,n$
and we only displayed columns labelled by $j$;
plus
\beql{explamt}
\sum_{j=i+1}^{n}\
\begin{vmatrix}
\ \ \dots&(\la_i-\la_1)\tss \Ec_{1\tss i}&\dots\ \ \\[0.2em]
\ \ \dots&\dots&\dots\ \  \\[0.2em]
\ \ \dots&(\la_{i}-\la_{i-1})\tss \Ec_{i-1\ts i} &\dots\ \  \\[0.2em]
\ \ \dots&0&\dots\ \ \\[0.2em]
\ \ \dots&\dots&\dots\ \  \\[0.2em]
\ \ \dots&0&\dots\ \
\end{vmatrix}
\quad+\quad
\sum_{j=i+1}^{n-1}(\la_j-\la_i)\ts \sum_{m=j+1}^n\
\begin{vmatrix}
\ \ \dots&0&\dots\ \ \\[0.2em]
\ \ \dots&\dots&\dots\ \  \\[0.2em]
\ \ \dots&\Ec_{i\tss m} &\dots\ \  \\[0.2em]
\ \ \dots&\dots&\dots\ \  \\[0.2em]
\ \ \dots&0&\dots\ \
\end{vmatrix},
\eeq
where the rows and columns of
the column-determinants are labelled by $1,\dots,\wh\jmath,\dots,n$
and we displayed columns $i$ and $m$, respectively; together with
\beql{explamtau}
\sum_{j=1}^{i-1}\ts\la_j\ts
\begin{vmatrix}
\ \ \dots&0&\dots\ \ \\[0.2em]
\ \ \dots&\dots&\dots\ \  \\[0.2em]
\ \ \dots&\tau_i-\tau_j &\dots\ \  \\[0.2em]
\ \ \dots&\dots&\dots\ \  \\[0.2em]
\ \ \dots&0&\dots\ \
\end{vmatrix}
\quad+\quad
\sum_{j=i+1}^{n}\ts\la_i\ts
\begin{vmatrix}
\ \ \dots&0&\dots\ \ \\[0.2em]
\ \ \dots&\dots&\dots\ \  \\[0.2em]
\ \ \dots&\tau_i-\tau_j &\dots\ \  \\[0.2em]
\ \ \dots&\dots&\dots\ \  \\[0.2em]
\ \ \dots&0&\dots\ \
\end{vmatrix},
\eeq
where the rows and columns of
the column-determinants are labelled by $1,\dots,\wh\jmath,\dots,n$
in the first sum and by $1,\dots,\wh\imath,\dots,n$ in the second, and
$\tau_i-\tau_j$ occurs as the $(i,i)$ and $(j,j)$-entry, respectively.
It is understood that the remaining non-displayed entries of all
column-determinants coincide with the respective entries of the matrix $\Ec$.

The final part of the arguments is a straightforward verification
that the contribution of
each of the above column-determinants to
the coefficient of $x^{n-k}$
in the polynomial $E^{(0)}_{i\tss i}[1]\tss\cdet\ts\Ec$
is a polynomial in $u$ whose degree is less than \eqref{lowb}.

\paragraph{Generators $E^{(0)}_{i\tss i}[s]$ for $s\geqslant 2${\rm .}}
Note that due to the fact that the element $E^{(0)}_{1\tss 1}+\dots+E^{(0)}_{n\tss n}$
belongs to the kernel of the bilinear form $\langle\ts\ts,\ts\rangle$, this
step is not necessary for the Takiff case, where $\la_1=\dots=\la_n$.
Indeed, all elements $E^{(0)}_{1\tss 1}[s]+\dots+E^{(0)}_{n\tss n}[s]$ are central
in the Lie algebra $\wh\agot$, whereas in $\agot$ we have
$E^{(0)}_{i\tss i}-E^{(0)}_{i+1\ts i+1}=[E^{(0)}_{i\ts i+1},E^{(0)}_{i+1\tss i}]$.
Therefore, the required vanishing properties of $\cdet\ts\Ec$ for the generators
$E^{(0)}_{i\tss i}[s]$ for $s\geqslant 2$ follow from the previously considered cases.

Returning to the general values of the $\la_i$, note the relations
\beql{rered}
s\ts E^{(p)}_{i\tss i}[s+1]=\big[\De,E^{(p)}_{i\tss i}[s]\big]
\eeq
implied by \eqref{bde}. Therefore the required properties
$E^{(0)}_{i\tss i}[s]\tss \phi^{(r)}_k=0$ for $s\geqslant 2$ follow from Lemma~\ref{lem:acde}
by a straightforward induction on $s$ and $k$, where the cases $s=1$ and $k=1$ are taken
as the induction bases.

\bpf[Proof of Lemma~\ref{lem:acde}]
The commutator of two derivations $\De$ and $T$ of $V(\agot)\cong\U\big(t^{-1}\agot[t^{-1}]\big)$
is another derivation
$
[\De,T]=2\tss d,
$
where the action of $d$ is determined by the relations
\ben
\big[d,X[r]\tss\big]=r\ts X[r],\qquad
X\in\agot.
\een
We also have the commutation relation $[d,T]=-T$.

Applying the operator $\De$ to the column-determinant
$\cdet\ts\Ec$ we get the expansion
\beql{deact}
\De \ts\cdet\ts\Ec=
\sum_{j=1}^n\ \begin{vmatrix}
\Ec_{1\tss 1}&\dots&-\Ec_{1\tss j}[0]&\dots&\Ec_{1\tss n}\\[0.2em]
\dots&\dots&\dots&\dots &\dots \\[0.2em]
\dots&\dots&2\la_j\tss d-\Ec_{j\tss j}[0] &\dots &\dots \\[0.2em]
\dots&\dots&\dots&\dots &\dots \\[0.2em]
\Ec_{n\tss 1}&\dots&-\Ec_{n\tss j}[0]&\dots&\Ec_{n\tss n}
\end{vmatrix},
\eeq
where we use notation \eqref{eijze}. The next step is to evaluate
these column-determinants in the vacuum module in the same way as for the
determinants \eqref{difftwo} above, with the use of the relations
$
[d,\ts\Ec_{ij}]=-\Ec^{\circ}_{ij}
$
with $\Ec^{\circ}_{ij}=\Ec_{ij}-\de_{ij}\tss x$.
The contribution arising from the action of $d$ equals
\ben
\sum_{j=1}^{n-1} 2\la_j\ts\sum_{m=j+1}^n\ \begin{vmatrix}
\Ec_{1\tss 1}&\dots&-\Ec^{\circ}_{1\tss m}&\dots&\Ec_{1\tss n}\\[0.2em]
\dots&\dots&\dots&\dots &\dots \\[0.2em]
\Ec_{n\tss 1}&\dots&-\Ec^{\circ}_{n\tss m}&\dots&\Ec_{n\tss n}
\end{vmatrix},
\een
where rows and columns of the column-determinants
are labelled by the numbers $1,\dots,\wh\jmath,\dots,n$,
and we displayed column $m$ in the middle.

To evaluate the column-determinants in \eqref{deact}
containing the terms $\Ec_{i\tss j}[0]$, apply formulas \eqref{eact}.
It is clear from these evaluations of the expression in \eqref{deact}
that the coefficient of $x^{n-k}$
is a polynomial in $u$ of degree not exceeding the number
given in \eqref{lowb}. This implies relations \eqref{dephize}
thus proving the first part of the lemma.
Furthermore, in proving its the second part, we only need to keep
the leading component in $u$ of degree equal to \eqref{lowb}.
Apply the above formulas obtained for the evaluation of
determinants \eqref{difftwo}, including suitable permutations of rows and columns
of column-determinants, which are valid
modulo lower degree terms in $u$. As a result, we can represent
the contribution arising from the terms $\Ec_{i\tss j}[0]$
as the sum of the expressions
\beql{laxle}
\sum_{j=1}^{n-1}\sum_{m=j+1}^n\
\begin{vmatrix}
\Ec_{1\tss 1}&\dots&\la_1\tss \Ec^{\circ}_{1\tss m}&\dots&\Ec_{1\tss n}\\[0.2em]
\dots&\dots&\dots&\dots &\dots \\[0.2em]
\dots&\dots&\la_{j-1}\tss \Ec^{\circ}_{j-1\ts m} &\dots &\dots \\[0.2em]
\dots&\dots&\la_j\tss \Ec^{\circ}_{j+1\ts m} &\dots &\dots \\[0.2em]
\dots&\dots&\dots&\dots &\dots \\[0.2em]
\Ec_{n\tss 1}&\dots&\la_j\tss \Ec^{\circ}_{n\tss m}&\dots&\Ec_{n\tss n}
\end{vmatrix},
\eeq
where the rows and columns of
the column-determinants are labelled by $1,\dots,\wh\jmath,\dots,n$;
plus the sum
\beql{laxra}
\sum_{j=1}^{n-1}\sum_{i=j+1}^n\
\begin{vmatrix}
\Ec_{1\tss 1}&\dots&\la_1\tss \Ec^{\circ}_{1\tss j}&\dots&\Ec_{1\tss n}\\[0.2em]
\dots&\dots&\dots&\dots &\dots \\[0.2em]
\dots&\dots&\la_{j-1}\tss \Ec^{\circ}_{j-1\ts j} &\dots &\dots \\[0.2em]
\dots&\dots&\la_j\tss \Ec^{\circ}_{j\tss j} &\dots &\dots \\[0.2em]
\dots&\dots&\dots&\dots &\dots \\[0.2em]
\Ec_{n\tss 1}&\dots&\la_j\tss \Ec^{\circ}_{n\tss j}&\dots&\Ec_{n\tss n}
\end{vmatrix},
\eeq
where the rows and columns of
the column-determinants are labelled by $1,\dots,\wh\imath,\dots,n$;
together with
\beql{taudi}
\sum_{j=1}^{n-1}\sum_{i=j+1}^n\
\begin{vmatrix}
\Ec_{1\tss 1}&\dots&0&\dots&\Ec_{1\tss n}\\[0.2em]
\dots&\dots&\dots&\dots &\dots \\[0.2em]
\dots&\dots&\la_j\tss (\tau_j-\tau_i) &\dots &\dots \\[0.2em]
\dots&\dots&\dots&\dots &\dots \\[0.2em]
\Ec_{n\tss 1}&\dots&0&\dots&\Ec_{n\tss n}
\end{vmatrix},
\eeq
where the rows and columns of
the column-determinants are labelled by $1,\dots,\wh\jmath,\dots,n$
and $\la_j\tss (\tau_j-\tau_i)$ occurs as the $(i,i)$-entry.

As with the determinants of the form \eqref{explam} and \eqref{explamtau}
considered above, we can see
that no contribution to
the leading component in $u$ of degree equal to \eqref{lowb}
comes from the expression \eqref{taudi}, while
modulo lower degree terms in $u$, both expressions \eqref{laxle} and \eqref{laxra}
can be simplified by replacing all coefficients $\la_1,\dots,\la_{j-1}$ with $\la_j$.

Thus, summarizing the arguments, we may conclude that
for $r=\la_{n-k+2}+\dots+\la_n-k+1$
the image
$\De(\phi^{(r)}_k)$ coincides with the coefficient of $u^r$
of the polynomial in $u$ which occurs as the coefficient of $x^{n-k}$
in the expression
\ben
{}-\sum_{j=1}^{n-1} \la_j\ts\sum_{m=j+1}^n\ \begin{vmatrix}
\Ec_{1\tss 1}&\dots&\Ec^{\circ}_{1\tss m}&\dots&\Ec_{1\tss n}\\[0.2em]
\dots&\dots&\dots&\dots &\dots \\[0.2em]
\Ec_{n\tss 1}&\dots&\Ec^{\circ}_{n\tss m}&\dots&\Ec_{n\tss n}
\end{vmatrix}-\sum_{j=1}^{n-1} \la_j\ts\sum_{i=j+1}^n\ \begin{vmatrix}
\Ec_{1\tss 1}&\dots&\Ec^{\circ}_{1\tss j}&\dots&\Ec_{1\tss n}\\[0.2em]
\dots&\dots&\dots&\dots &\dots \\[0.2em]
\Ec_{n\tss 1}&\dots&\Ec^{\circ}_{n\tss j}&\dots&\Ec_{n\tss n}
\end{vmatrix},
\een
where the rows and columns of
the column-determinants are labelled by $1,\dots,\wh\jmath,\dots,n$
in the first sum and by $1,\dots,\wh\imath,\dots,n$ in the second.

The coefficient of $x^{n-k}$ in the above expression is a linear combination
of the principal minors $\Ec^{a_1,\dots,a_{k-1}}_{}$ of the matrix $\Ec$
corresponding to the rows and columns $a_1<\dots<a_{k-1}$. In particular,
observe that the linear combination contains the principal minor $\Ec^{n-k+2,\dots,n}_{}$
with the coefficient
\beql{coek}
-(k-1)\tss(\la_1+\dots+\la_{n-k+1}).
\eeq
To apply
the formula for $\phi^{(r)}_{k-1}$ provided by \eqref{polt}
and complete the proof of the lemma we need to verify that
all other principal minors containing $u^r$ occur in the
linear combination with the same coefficient.

Define the parameters $1\leqslant s\leqslant n-k+2$ and $n-k+2\leqslant p\leqslant n$
by the conditions
\ben
\la_1\leqslant\dots\leqslant\la_{s-1}<\la_{s}=\dots=\la_{n-k+2}=\dots=\la_p
<\la_{p+1}\leqslant\dots\leqslant \la_n.
\een
The principal minors containing $u^r$ must have the form
$\Ec^{a_1,\dots,a_{q},p+1,\dots,n}_{}$ for some
\ben
s\leqslant a_1<\dots<a_q\leqslant p,\qquad q=p-n+k-1.
\een
The coefficient of such a minor in the linear combination equals
\ben
-\sum_{i=1}^{q+1}\ts(k-i)(\la^{}_{a_{i-1}+1}+\dots+\la^{}_{a_i-1})
-\sum_{i=1}^q\ts(p-q-a_i+i)\la^{}_{a_i},
\een
where we set $a_0=0$ and $a_{q+1}=p+1$. This coefficient coincides with
\eqref{coek} thus completing the proof of the lemma.
\epf

To establish the required vanishing properties of the elements $\phi^{(r)}_k$
with respect to the remaining family of generators $E^{(p)}_{i\tss i}[s]$
with $s\geqslant 0$, observe that
in view of Lemma~\ref{lem:acde} and relations \eqref{rered}, it is sufficient
to consider the values $s=0$ and $s=1$.

\paragraph{Generators $E^{(p)}_{i\tss i}[0]$ with $p=0,1,\dots,\la_i-1${\rm .}}

For a fixed value or $p$ we have the commutation relations
\ben
u^p\tss \big[E^{(p)}_{i\tss i}[0],\Ec_{m\tss l}\big]=
\de_{m\tss i}\tss (\Ec_{i\tss l}-\wt \Ec_{i\tss l})
-\de_{i\tss l}\tss (\Ec_{m\tss i}-\wh \Ec_{m\tss i}),
\een
where we set $\wt\Ec_{i\tss i}=\wh \Ec_{i\tss i}=0$,
\ben
\wt\Ec_{i\tss l}=\begin{cases}
E^{(\la_l-\la_i)}_{i\tss l}[-1]\tss u^{\la_l-\la_i}
+\dots+E^{(\la_l-\la_i+p-1)}_{i\tss l}[-1]\tss u^{\la_l-\la_i+p-1}\ \ &\text{for}\quad i<l,\\[0.3em]
E^{(0)}_{i\tss l}[-1]
+\dots+E^{(p-1)}_{i\tss l}[-1]\tss u^{p-1}\quad&\text{for}\quad i>l,\ \ p<\la_l,\\[0.3em]
E^{(0)}_{i\tss l}[-1]
+\dots+E^{(\la_l-1)}_{i\tss l}[-1]\tss u^{\la_l-1}\quad&\text{for}
\quad i>l,\ \ p\geqslant\la_l
\end{cases}
\een
and
\ben
\wh\Ec_{m\tss i}=\begin{cases}
E^{(\la_i-\la_m)}_{m\tss i}[-1]\tss u^{\la_i-\la_m}
+\dots+E^{(\la_i-\la_m+p-1)}_{m\tss i}[-1]\tss u^{\la_i-\la_m+p-1}
\ \ &\text{for}\quad m<i,\ \   p<\la_m,\\[0.3em]
E^{(\la_i-\la_m)}_{m\tss i}[-1]\tss u^{\la_i-\la_m}
+\dots+E^{(\la_i-1)}_{m\tss i}[-1]\tss u^{\la_i-1}
\ \ &\text{for}\quad m<i,\ \   p\geqslant\la_m,\\[0.3em]
E^{(0)}_{m\tss i}[-1]
+\dots+E^{(p-1)}_{m\tss i}[-1]\tss u^{p-1}
\ \ &\text{for}\quad m>i.
\end{cases}
\een
Hence the element $E^{(p)}_{i\tss i}[0]\ts\cdet\ts\Ec$
equals the difference of two column-determinants
\ben
\begin{vmatrix}
\Ec_{1\ts 1}&\dots&\Ec_{1\ts n}\\[0.2em]
\dots&\dots &\dots \\[0.2em]
\Ec_{i\tss 1}-\wt\Ec_{i\tss 1}&\dots&\Ec_{i\tss n}-\wt\Ec_{i\tss n}\\[0.2em]
\dots&\dots&\dots \\[0.2em]
\Ec_{n\ts 1}&\dots&\Ec_{n\ts n}
                \end{vmatrix}
\quad - \quad
                \begin{vmatrix}
\Ec_{1\ts 1}&\dots&\Ec_{1\ts i}-\wh\Ec_{1\ts i}&\dots&\Ec_{1\ts n}\\[0.2em]
\dots&\dots&\dots&\dots  &\dots \\[0.2em]
\Ec_{n\ts 1}&\dots&\Ec_{n\ts i}-\wh\Ec_{n\ts i}&\dots&\Ec_{n\ts n}
                \end{vmatrix},
\een
which simplifies to
\beql{dirze}
\begin{vmatrix}
\Ec_{1\ts 1}&\dots&\wh\Ec_{1\ts i}&\dots&\Ec_{1\ts n}\\[0.2em]
\dots&\dots&\dots&\dots  &\dots \\[0.2em]
\Ec_{n\ts 1}&\dots&\wh\Ec_{n\ts i}&\dots&\Ec_{n\ts n}
                \end{vmatrix}
\quad - \quad
\begin{vmatrix}
\Ec_{1\ts 1}&\dots&\Ec_{1\ts n}\\[0.2em]
\dots&\dots &\dots \\[0.2em]
\wt\Ec_{i\tss 1}&\dots&\wt\Ec_{i\tss n}\\[0.2em]
\dots&\dots&\dots \\[0.2em]
\Ec_{n\ts 1}&\dots&\Ec_{n\ts n}
                \end{vmatrix}.
\eeq
It remains to verify that each of these column-determinants
is a polynomial in $x$ with the property
that the coefficient of $x^{n-k}$ is a polynomial in $u$ whose degree
is less that $p$ plus the expression in \eqref{lowb};
the calculations are quite similar to the previous cases.

\paragraph{Generators $E^{(p)}_{i\tss i}[1]$ with $p=1,\dots,\la_i-1${\rm .}}

We proceed as in the case $p=0$ considered above and use similar calculations
with some simplifications due to the fact that they do not involve
values of the bilinear form.
We have
the relations
\ben
u^p\tss \big[E^{(p)}_{i\tss i}[1], \Ec_{m\tss l}\big]=
\de_{m\tss l}\tss \la_m\tss u^p\tss E^{(p)}_{i\tss i}[0]
+\de_{m\tss i}\tss \ol\Ec_{i\tss l}[0]-\de_{i\tss l}\tss \ol\Ec_{m\tss i}[0],
\een
where we use the notation
\ben
\ol\Ec_{i\tss j}[0]=\begin{cases}
E^{(\la_j-\la_i+p)}_{i\tss j}[0]\tss u^{\la_j-\la_i+p}
+\dots+E^{(\la_j-1)}_{i\tss j}[0]\tss u^{\la_j-1}\ \ &\text{for}\quad i<j,\\[0.3em]
E^{(p)}_{i\tss j}[0]\tss u^{p}
+\dots+E^{(\la_j-1)}_{i\tss j}[0]\tss u^{\la_j-1}\ \ &\text{for}\quad i>j,
\end{cases}
\een
where empty sums are understood as being equal to zero. Relations \eqref{eoneact}
hold in the same form with $E^{(0)}_{i\tss i}[1]$
replaced by $u^p\tss E^{(p)}_{i\tss i}[1]$ so that the commutator
$\big[u^p\tss E^{(p)}_{i\tss i}[1],\cdet\ts\Ec\big]$ equals the sum
\ben
\sum_{j=1}^n\ts
\begin{vmatrix}
\Ec_{1\tss 1}&\dots&0&\dots&\Ec_{1\tss n}\\[0.2em]
\dots&\dots&\dots&\dots &\dots \\[0.2em]
\dots&\dots&\la_j\ts u^p\tss E^{(p)}_{i\tss i}[0] &\dots &\dots \\[0.2em]
\dots&\dots&\dots&\dots &\dots \\[0.2em]
\Ec_{n\tss 1}&\dots&0&\dots&\Ec_{n\tss n}
                \end{vmatrix}
\een
plus the difference of two column-determinants
\ben
\begin{vmatrix}
\Ec_{1\tss 1}&\dots&\Ec_{1\tss n}\\[0.2em]
\dots&\dots&\dots \\[0.2em]
\ol\Ec_{i\tss 1}[0]&\dots&\ol\Ec_{i\tss n}[0] \\[0.2em]
\dots&\dots&\dots \\[0.2em]
\Ec_{n\tss 1}&\dots&\Ec_{n\tss n}
                \end{vmatrix}
\quad - \quad
                \begin{vmatrix}
\Ec_{1\tss 1}&\dots&\ol\Ec_{1\ts i}[0]&\dots&\Ec_{1\tss n}\\[0.2em]
\dots&\dots&\dots&\dots &\dots \\[0.2em]
\Ec_{n\tss 1}&\dots&\ol\Ec_{n\ts i}[0]&\dots&\Ec_{n\tss n}
                \end{vmatrix}.
\een
It turns out that in contrast to the case $p=0$,
each of these column-determinants
is a polynomial in $x$ with the property
that the coefficient of $x^{n-k}$ is a polynomial in $u$ whose degree
is less that $p$ plus the expression in \eqref{lowb}.
Indeed, this follows easily from the respective counterparts
of the decompositions \eqref{midt} and \eqref{eact},
because the resulting column-determinants retain the property that the degree of $u$
of any entry in column $l$ does not exceed $\la_l-1$.

We have thus proved that
the coefficients $\phi^{(r)}_k$ with $k=1,\dots,n$ satisfying conditions
\eqref{conda} belong to the center $\z(\wh\agot)$
of the vertex algebra $V(\agot)$.
To prove the second part of the theorem,
equip the vector space $V(\agot)\cong \U(t^{-1}\agot[t^{-1}])$
with the canonical filtration of the universal enveloping algebra and
consider the
symbols $\bar\phi^{(r)}_k$ of the coefficients
in the associated graded algebra $\Sr(t^{-1}\agot[t^{-1}])=\gr \U(t^{-1}\agot[t^{-1}])$.
The symbols belong to the subalgebra of $\agot[t]$-invariants in $\Sr(t^{-1}\agot[t^{-1}])$,
and it is clear from the definition of the elements $\phi^{(r)}_k\in V(\agot)$ that
each $\bar\phi^{(r)}_k$ coincides with the image of
a certain element of the algebra of $\agot$-invariants
$\Sr(\agot)^{\agot}$ under the embedding $\imath:X\mapsto X[-1]$ for $X\in\agot$.
More precisely, introduce polynomials in $u$ with coefficients
in $\Sr(\agot)$ by
\ben
\overline\Ec_{ij}(u)=\begin{cases}E^{(0)}_{ij}+\dots+E^{(\la_j-1)}_{ij}\ts u^{\la_j-1}
&\text{if}\quad i\geqslant j,\\[0.4em]
E^{(\la_j-\la_i)}_{ij}\tss u^{\la_j-\la_i}+\dots+E^{(\la_j-1)}_{ij}\ts u^{\la_j-1}
&\text{if}\quad i< j
\end{cases}
\een
and write the determinant
\ben
\det\left[\begin{matrix}
x+\overline\Ec_{11}(u)&\overline\Ec_{12}(u)&\dots&\overline\Ec_{1n}(u)\\[0.4em]
\overline\Ec_{21}(u)&x+\overline\Ec_{22}(u)&\dots&\overline\Ec_{2n}(u)\\
\vdots&\vdots& \ddots&\vdots     \\
\overline\Ec_{n1}(u)&\overline\Ec_{n2}(u)&\dots&x+\overline\Ec_{n\tss n}(u)
                \end{matrix}\right]
\een
as a polynomial in $x$,
\ben
x^n+\overline\Phi_1(u)\tss x^{n-1}+\dots+\overline\Phi_n(u)\qquad\text{with}\quad
\overline\Phi_k(u)=\sum_r\overline\Phi^{(r)}_k\ts u^r.
\een
Then under the embedding $\imath$ we have $\imath:\overline\Phi^{(r)}_k\mapsto \bar\phi^{(r)}_k$.
By \cite[Thm~4.1]{bb:ei}, the
coefficients $\overline\Phi^{(r)}_k$ with $k=1,\dots,n$ and
$r$ satisfying conditions \eqref{conda}
are algebraically independent generators of the algebra $\Sr(\agot)^{\agot}$.
The same generators were produced
in an earlier work \cite[Conj.~4.1]{ppy:si} as a conjecture, which was later confirmed in
\cite[Thm~7]{o:sp}. In the terminology of \cite{ppy:si},
these generators arise from the {\em good generating system} of $\Sr(\gl_N)^{\gl_N}$
for the nilpotent element $e$, formed by the coefficients of the
characteristic polynomial of the matrix $[E_{ij}]_{i,j=1}^N$.
This generating system satisfies
the assumptions of \cite[Thm~3.2]{ap:qm} which implies that
$\Sr(t^{-1}\agot[t^{-1}])^{\agot[t]}$ is a polynomial algebra
in the variables $T^j \bar\phi^{(r)}_k$ with $j\geqslant 0$,
where $k=1,\dots,n$ and $r$ satisfies conditions \eqref{conda}.
This completes the proof of Theorem~\ref{thm:suga}, since under these conditions
on the indices, the elements $T^j \phi^{(r)}_k$ are algebraically
independent generators of the algebra $\z(\wh\agot)$.

\section*{Acknowledgments}

The support of
the Australian Research Council, grant DP180101825
is gratefully acknowledged.

%\newpage
\bigskip\bigskip

\small

\noindent
School of Mathematics and Statistics\newline
University of Sydney,
NSW 2006, Australia\newline
alexander.molev@sydney.edu.au

\end{document}